\documentclass{CSML}
\pdfoutput=1

\usepackage{lastpage}
\newcommand{\dOi}{13(4:18)2017}
\lmcsheading{}{1--\pageref{LastPage}}{}{}%
{Sep.~06, 2016}{Nov.~27, 2017}{}

\pdfoutput=1 

\usepackage{hyperref}
\hypersetup{hidelinks}
\usepackage{latexsym,amsmath,amssymb,stmaryrd}
\usepackage{color}
\usepackage[all]{xy}

\theoremstyle{plain}

\newcommand\pow{\mathbb{P}}
\newcommand{\real}{\mathbb{R}}
\newcommand{\creal}{\overline{\mathbb{R}}_+}
\newcommand{\Rp}{\mathbb{R}_+}
\newcommand{\Rl}{\mathbb{R}_\ell}
\DeclareMathOperator{\upc}{\uparrow\!}
\DeclareMathOperator{\dc}{\downarrow\!}
\newcommand\nat{\mathbb{N}}
\newcommand\dreal{\text{\textsf{\upshape d}}_\real}
\newcommand\dRl{\text{\textsf{\upshape d}}_\ell}
\newcommand\Idl{\mathbf{I}}
\newcommand\limp{\mathrel{\Rightarrow}}

\newcommand\FB{\mathbf S}
\newcommand\RIdl{\mathbf{RI}}
\newcommand\uuarrow{\rlap{$\uparrow$}\raise.5ex\hbox{$\uparrow$}}
\newcommand\ddarrow{\rlap{$\downarrow$}\raise.5ex\hbox{$\downarrow$}}
\newcommand\Dc{\mathop{\Downarrow}}
\newcommand\identity[1]{\mathrm{id}_{#1}}

\newcommand\Hoare{{\mathcal H}}

\newcommand\Plotkin{\mathcal P\ell}
\newcommand\Plotkinn{\Plotkin_{\mathcal V}}
\newcommand\PV\Plotkinn 

\newcommand{\QMet}{\mathbf{QMet}}
\newcommand{\Dcpo}{\mathbf{Dcpo}}
\newcommand{\WC}{\mathbf{Y}}
\newcommand{\Ord}{\mathbf{Ord}}

\begin{document}

\title[Formal Balls]{A Few Notes on Formal Balls}

\author[Goubault-Larrecq]{Jean Goubault-Larrecq}	
\address{LSV, ENS Paris-Saclay, CNRS, Universit\'e Paris-Saclay, 94230 Cachan, France}	
\email{goubault@lsv.fr}  

\author[Ng]{Kok Min Ng}	
\address{National Institute of Education, Nanyang Technological
  University, Singapore}	
\email{ngkokmin@gmail.com}  



\keywords{Quasi-metric, formal balls}
\subjclass{Mathematics of computing - continuous mathematics -
  topology - point-set topology}


\begin{abstract}
  \noindent Using the notion of formal ball, we present a few new
  results in the theory of quasi-metric spaces.  With no specific
  order: every continuous Yoneda-complete quasi-metric space is sober
  and convergence Choquet-complete hence Baire in its $d$-Scott
  topology; for standard quasi-metric spaces, algebraicity is
  equivalent to having enough center points; on a standard
  quasi-metric space, every lower semicontinuous $\creal$-valued
  function is the supremum of a chain of Lipschitz Yoneda-continuous
  maps; the continuous Yoneda-complete quasi-metric spaces are exactly
  the retracts of algebraic Yoneda-complete quasi-metric spaces; every
  continuous Yoneda-complete quasi-metric space has a so-called
  quasi-ideal model, generalizing a construction due to K. Martin.
  The point is that all those results reduce to domain-theoretic
  constructions on posets of formal balls.
\end{abstract}

\maketitle

\section*{Introduction}
\label{sec:into}

In his gem of a paper on how to write Mathematics
\cite[Section~2]{Halmos:math}, Paul Halmos recommends to ``say
something''.  He then comments on books and papers that violate this
principle by either saying nothing or saying too many things.  The
present paper may appear to say too many, relatively random, things.
On the contrary, let us stress that the unique idea of the present
paper can be summarized by the motto: ``Formal balls are the essence
of quasi-metric spaces''.

We will explain all terms in Section~\ref{sec:basics-quasi-metric}.

The first author has been convinced of the truth of that motto while
writing the book \cite{JGL-topology}, and most of its Chapter~7 arises
from that conviction.  Several papers had already been based on that
premise
\cite{Vickers:completion:gms:I,Vickers:completion:gms:II,AHPR:formal:ball,RV:formal:ball,KW:formal:ball}.
The first author has given a few talks on the topic, in particular at
the Domains XII conference in Cork, Ireland, 2015.  However, most of
it has already been published, and only a few crumbs remain to offer
the reader.  We hope those are interesting crumbs.

We outline a notion of \emph{standard} quasi-metric space in
Section~\ref{sec:stand-quasi-metr}.  This is a natural notion, and all
Yoneda-complete quasi-metric spaces, as well as all metric spaces, are
standard.  We give a simple characterization of Waszkiewicz's
\emph{continuous} Yoneda-complete quasi-metric spaces in
Section~\ref{sec:cont-yoneda-compl}, and this suggests a definition of
continuous, not necessarily complete, quasi-metric spaces: namely, the
standard quasi-metric spaces whose poset of formal balls is
continuous.  All metric spaces are continuous in this sense, in
particular.

The above crumbs can be considered as additional basic facts on
quasi-metric spaces, complementing
Section~\ref{sec:basics-quasi-metric}.  Those facts have not appeared
earlier in the literature, as far as we know.  In
Section~\ref{sec:choquet-completeness}, we grab low-hanging fruit and
show that every continuous Yoneda-complete quasi-metric space is
sober, and also convergence Choquet-complete hence Baire.  In
Section~\ref{sec:algebr-yoneda}, we characterize those standard
quasi-metric spaces that are \emph{algebraic}, as those that have
\emph{enough center points}.  The latter is a simple condition on the
Scott topology on the poset of formal balls.  In
Section~\ref{sec:cont-lipsch-real}, we look at morphisms, and first
show that every lower semicontinuous map from a standard quasi-metric
space to $\creal$ is a pointwise supremum of a chain of Lipschitz
Yoneda-continuous maps.  This generalizes a standard construction on
metric spaces, and involves defining the appropriate variant of the
distance of a point to a closed set.  Then, in
Section~\ref{sec:cont-algebr-quasi}, we show that the continuous
Yoneda-complete quasi-metric spaces are exactly the retracts of
algebraic Yoneda-complete quasi-metric spaces, generalizing a similar
result in the theory of dcpos.

A final crumb, in Section~\ref{sec:quasi-ideal-models}, explores the
notion of quasi-ideal domains: algebraic dcpos whose finite elements
are below all non-finite elements.  Using a variant of a construction
due to K. Martin, we show that every continuous Yoneda-complete
quasi-metric space has a quasi-ideal model; and that the spaces that
have an $\omega$-quasi-ideal model are exactly M. de Brecht's
quasi-Polish spaces.

We conclude in Section~\ref{sec:concl-open-probl}.

\section{Basics on Quasi-Metric Spaces}
\label{sec:basics-quasi-metric}

Let $\creal$ be the set of extended non-negative reals.  A
\emph{quasi-metric} on a set $X$ is a map
$d \colon X \times X \to \creal$ satisfying: $d (x, x)=0$;
$d (x, z) \leq d (x, y) + d (y, z)$ (triangular inequality);
$d(x,y)=d(y,x)=0$ implies $x=y$.

The relation $\leq^d$ defined by $x \leq^d y$ if and only if
$d (x,y)=0$ is then an ordering, which turns out to be the
specialization ordering of the open ball topology.  The latter is
generated by the open balls $B^d_{x, <r}$, $x \in X$, $r > 0$, as
basic open sets.

\begin{exa}
  \label{ex:dreal}
  On any subset of $\real \cup \{+\infty\}$, define the quasi-metric
  $\dreal$ by $\dreal (x, y) = 0$ if $x \leq y$,
  $\dreal (+\infty, y) = +\infty$ if $y \neq +\infty$,
  $\dreal (x, y) = x-y$ if $x > y$ and $x \neq +\infty$.  This is a
  quasi-metric with specialization ordering $\leq^{\dreal}$ equal to
  the usual ordering $\leq$.  \qed
\end{exa}

A \emph{formal ball} is a pair $(x, r)$ where $x \in X$ and
$r \in \Rp$.  This is just syntax for an actual ball: $x$ is the
\emph{center}, and $r$ is the \emph{radius}.  Formal balls are ordered
by $(x, r) \leq^{d^+} (y, s)$ if and only if $d (x, y) \leq r-s$.
Note that this implies $r \geq s$, in particular.

The poset of formal balls $\mathbf B (X, d)$ has many serendipitous
properties.  Most of them are described and proved in
\cite{JGL-topology}, and we will recapitulate the most fundamental ones
here.  In the remaining sections, we will proceed to state a few
new results that stem from the study of $\mathbf B (X, d)$.

Although the book \cite{JGL-topology} is a good source of information,
we would not like to insinuate that the first author of the present
paper is the author of the theory of formal balls.  Formal balls were
introduced by Klaus Weihrauch and Ulrich Schreiber
\cite{WS:formal:ball}.  Reinhold Heckmann and Abbas Edalat showed why
they were so important in the metric case \cite{EH:comp:metric}.  In
the general, quasi-metric case, we would like to stress the import of
Mateusz Kostanek and Pawe\l{} Waszkiewicz \cite{KW:formal:ball}, who
showed that $X, d$ is Yoneda-complete if and only if
$\mathbf B (X, d)$ is a dcpo ---~such an important characterization
that we will actually take it as a \emph{definition} of
Yoneda-completeness, although we will recapitulate the standard
definition in Section~\ref{sec:algebr-yoneda}.  (Kostanek and
Waszkiewicz consider the more general case of quasi-metrics with
values in a $\mathcal Q$-category instead of $\creal$, a natural
extension if we start from Lawvere's view of metric spaces as enriched
categories \cite{Lawvere:metric}.)  The supremum of a directed family
of formal balls ${(x_i, r_i)}_{i \in I}$ is then equal to $(x, r)$
where $x$ is the so-called \emph{$d$-limit} of the Cauchy net
${(x_i)}_{i \in I, \leq}$---here $i \leq j$ iff
$(x_i, r_i) \leq^{d^+} (x_j, r_j)$.  We will again define what a
$d$-limit is, and what Cauchy nets are, in
Section~\ref{sec:algebr-yoneda}.  The point is that the $d$-limit is
independent of the radii $r_i$.  Also, we have the innocuous-looking
equality $r = \inf_{i \in I} r_i$, which will be a crucial property of
standard quasi-metric spaces (Section~\ref{sec:stand-quasi-metr}).

Salvador Romaguera and Oscar Valero proved an important theorem of the
same kind, stating that $X, d$ is Smyth-complete if and only if
$\mathbf B (X, d)$ is a continuous dcpo and its way-below relation is
the relation $\prec$, defined by $(x, r) \prec (y, s)$ if and only if
$d (x, y) < r-s$ \cite{RV:formal:ball}.  Smyth-completeness is a
stronger property than Yoneda-completeness, and means that every
(forward) Cauchy net in $X, d$ has a limit in the metric space $X,
d^{sym}$ (taking $d^{sym} (x,y) = \max (d (x, y), d (y, x))$).

It turns out that $\mathbf B (X, d)$ can also be given a natural
quasi-metric $d^+$, defined by
$d^+ ((x, r), (y, s)) = \max (d (x, y)-r+s, 0)$.  The associated
ordering is the ordering $\leq^{d^+}$ we have already defined:
$(x, r) \leq^{d^+} (y, s)$ if and only if $d (x, y) \leq r-s$.

With the open ball topology of $d^+$, $\mathbf B (X, d)$ is a
\emph{C-space}, as defined by Marcel Ern\'e \cite{Erne:ABC}.  In
general, a C-space is a $T_0$ topological space $Y$ that has the
following strong form of local compactness: for every $y \in Y$ and
every open neighborhood $V$ of $y$, there is a point $z \in V$ such
that $y$ is in the interior of the compact set $\upc z$.  (We write
$\upc z$ for the upward closure of $z$.)  This is particularly easy in
the case of $\mathbf B (X, d)$: for every open neighborhood $U$ (in
the open ball topology of $d^+$) of a formal ball $(x, r)$, $U$
contains an open ball $B^{d^+}_{(x, r), < \epsilon}$, and we can take
$z = (x, r+\epsilon/2)$; $(x, r)$ is in the interior of $\upc z$
since $(x, r)$ is in the open ball $B^{d^+}_{(x, r), < \epsilon/2}$,
which is included in $\upc z$.

Writing $z \prec y$ for ``$y$ is in the interior of $\upc z$'', that
makes every C-space an abstract basis, and conversely, any abstract
basis yields a C-space in a canonical way.  A sober C-space is exactly
the same as a continuous dcpo, the topology then being the Scott
topology, and $\prec$ being the way-below relation $\ll$.  In the case
of $\mathbf B (X, d)$, $\prec$ is the relation we have already
defined: $(x, r) \prec (y, s)$ if and only if $d (x, y) < r-s$.  It
follows that the Romaguera-Valero can be restated in the following
synthetic form: $X, d$ is Smyth-complete if and only if
$\mathbf B (X, d)$ is \emph{sober} in the open ball topology of $d^+$
\cite[Theorem~8.3.40]{JGL-topology}.  Sobriety is a fundamental
notion, and we refer the reader to \cite[Chapter~8]{JGL-topology} for
more information on the topic.

\begin{exa}
  \label{ex:Smyth-complete}
  With the quasi-metric $\dreal$, every closed interval $[a, b]$ is
  Smyth-complete.  In fact, $[a, b]$ is even symcompact, namely,
  compact in the open ball topology of the symmetrized metric
  $\dreal^{sym}$, and every symcompact space is Smyth-complete
  \cite[Lemma~7.2.21]{JGL-topology}.  \qed
\end{exa}

\begin{exa}
  \label{ex:Yoneda-complete:not}
  $\real, \dreal$ is not Smyth-complete, not even Yoneda-complete, and
  similarly for $\Rp, \dreal$.  Indeed, the chain of formal balls
  $(n, 0)$, $n \in \nat$, does not have a supremum.  \qed
\end{exa}

\begin{exa}
  \label{ex:Yoneda-complete:R}
  $\creal, \dreal$, on the contrary, is Yoneda-complete.  Noting that
  for $(x, r), (y, s) \in \mathbf (\creal, \dreal)$,
  $(x, r) \leq^{\dreal^+} (y, s)$ if and only if
  $r \geq s$ and $x-r \leq y-s$, the map $(x, r) \mapsto (x-r, -r)$
  defines an order isomorphism from $\mathbf B (\creal, \dreal)$ onto
  $C = \{(a, b) \in (\real \cup \{+\infty\}) \times (-\infty, 0] \mid
  a-b \geq 0\}$.
  Since $C$ is a Scott-closed subset of a continuous dcpo, it is
  itself a continuous dcpo.  We shall see in Example~\ref{ex:cont:R}
  that this can be used to immediately conclude that $\real, \dreal$
  is a so-called \emph{continuous} Yoneda-complete quasi-metric space.
  We can in fact say something stronger: $\creal, \dreal$ is algebraic
  \cite[Exercise~7.4.64]{JGL-topology}.  Since $(a, b) \ll (a', b')$
  in $C$ if and only if $a < a'$ and $b < b'$, we obtain that
  $(x, r) \ll (y, s)$ in $\mathbf B (\creal, \dreal)$ if and only if
  $x-r < y-s$ and $r > s$, if and only if $x \neq +\infty$ and
  $(x, r) \prec (y, s)$.  Hence the relations $\ll$ and $\prec$ differ
  (although by a narrow margin), showing that $\creal, \dreal$ is not
  Smyth-complete.  \qed
\end{exa}

\begin{exa}
  \label{ex:Yoneda-complete:Rl}
  The \emph{Sorgenfrey line} $\Rl$ is $\real, \dRl$, where
  $\dRl (x, y)$ is equal to $+\infty$ if $x > y$, and to $y-x$ if
  $x \leq y$.  Its specialization ordering $\leq^{\dRl}$ is equality.
  Its open ball topology is the topology generated by the half-open
  intervals $[x, x+r)$, and is a well-known counterexample in
  topology.  On formal balls, $(x, r) \leq^{\dRl^+} (y, s)$ if and
  only if $x \leq y$ and $x+r \geq y+s$.
  The map $(x, r) \mapsto (-x-r, x)$ therefore defines an order
  isomorphism from $\mathbf B (\real, \dRl)$ to
  $C_\ell = \{(a, b) \in \real^2 \mid a+b \leq 0\}$.  Although
  $\real^2$ is not a dcpo, $C_\ell$ is, as one can see by realizing
  that $C_\ell$ is the Scott-closed subset of the continuous dcpo
  $(\real \cup \{+\infty\})^2$ consisting of the pairs $(a, b)$ such
  that $a+b \leq 0$.  As such, $C$ is even a continuous dcpo, with
  $(a, b) \ll (a', b')$ if and only if $a<a'$ and $b<b'$.  The
  way-below relation on $\mathbf B (\real, \dRl)$ is given by
  $(x, r) \ll (y, s)$ if and only if $x+r > y+s$ and $x < y$, if and
  only if $(x, r) \prec (y, s)$ and $x < y$.  Since $\ll$ and $\prec$
  differ, $\Rl$ is not Smyth-complete.  One can show that it is even a
  non-algebraic Yoneda-complete quasi-metric space, see
  \cite{KW:formal:ball} or \cite[Exercise~7.4.73]{JGL-topology}.  This
  example is due to Kostanek and Waszkiewicz, who also show that it is
  continuous Yoneda-complete, a fact we shall retrieve in
  Example~\ref{ex:cont:Rl}.  \qed
\end{exa}

\begin{exa}
  \label{ex:<=}
  Every poset $X$ can be seen as a quasi-metric space by letting
  $d_\leq (x,y)=0$ if $x\leq y$, $d (x,y)=+\infty$ otherwise.  On
  formal balls, $(x, r) \leq^{d_\leq^+} (y, s)$ if and only if
  $x \leq y$ and $r \geq s$, so $(x, r) \mapsto (x, -r)$ defines an
  order isomorphism of $\mathbf B (X, d_\leq)$ onto
  $X \times (-\infty, 0]$.  It follows that $X$ is Yoneda-complete, as
  a quasi-metric space, if and only if $X$ is a dcpo.


  Note that $(x, r) \ll (y, s)$ in $\mathbf B (X, d)$ if and only if
  $x \ll y$ in $X$ and $r > s$, if and only if $(x, r) \prec (y, s)$
  and $x \ll y$.  In particular, $\ll$ and $\prec$ only coincide on
  formal balls when $\leq$ and $\ll$ coincide on $X$, that is, $X$ is
  Smyth-complete as a quasi-metric space if and only if it is an
  algebraic dcpo whose elements are all finite.  Such posets are
  exactly those that have the \emph{ascending chain condition}: every
  chain $x_1 \leq x_2 \leq \cdots \leq x_n \leq \cdots$ is finite.  \qed
\end{exa}

It has been argued that the proper topology one should take on a
quasi-metric space $X, d$ is not its open ball topology, but its
\emph{generalized Scott topology} \cite{BvBR:gms}.  We will use
another one, which is arguably simpler to understand, and coincides
with the latter in many cases.

The map $x \mapsto (x, 0)$ is an order embedding of $X, \leq^{d^+}$
into its poset of formal balls $\mathbf B (X, d)$.  Accordingly, we
shall consider $X$ as a subset of $\mathbf B (X, d)$.  The latter has
a natural topology, the Scott topology, and this induces a topology on
its subspace $X$:
\begin{defi}[$d$-Scott Topology]
  \label{defn:dScott}
  Let $X, d$ be a quasi-metric space.  The \emph{$d$-Scott topology}
  on $X$ is the induced topology from $\mathbf B (X, d)$, namely: the
  $d$-Scott opens subsets $U$ of $X$ are those such that there is a
  Scott-open subset $V$ of $\mathbf B (X, d)$ such that
  $U = \{x \mid (x, 0) \in V\}$.
\end{defi}
In general, the generalized Scott topology is finer than the $d$-Scott
topology on $X$ \cite[Exercise~7.4.51]{JGL-topology}.  The two
topologies coincide when $X$ is algebraic Yoneda-complete
\cite[Exercise~7.4.69]{JGL-topology}.

The reader might be worried at this point that we are relinquishing
the open ball topology in favor of a more exotic topology.  However,
the $d$-Scott topology is the \emph{same} as the open ball topology if
$X, d$ is \emph{metric} \cite[Proposition~7.4.46]{JGL-topology}.  The
two topologies also coincide if $X, d$ is Smyth-complete
\cite[Proposition~7.4.47]{JGL-topology}, and in fact as soon as all
the points of $X$ are $d$-finite (comment after Proposition~7.4.48,
op. cit.).

\begin{exa}
  \label{ex:Scott=dScott}
  For a poset $X$, seen as a quasi-metric space with the quasi-metric
  $d_\leq$ (Example~\ref{ex:<=}), the $d_\leq$-Scott topology is just
  the same as the Scott topology.  To show this, use the order
  isomorphism of $\mathbf B (X, d_\leq)$ onto $X \times (-\infty, 0]$.
  The $d$-Scott topology on $X$ is then the topology induced by the
  embedding $x \mapsto (x, 0)$ of $X$ into $X \times (-\infty, 0]$.
  We consider $X$ as a subset of $X \times (-\infty, 0]$ through that
  embedding.  For a Scott-open subset $U$ of $X$,
  $U \times (-\infty, 0]$ is Scott-open in $X \times (-\infty, 0]$,
  and its intersection with $X$ is $U$, so $U$ is $d$-Scott open.
  Conversely, any $d$-Scott open $U$ is of the form $V \cap X$ for
  some Scott open subset $V$ of $X \times (-\infty, 0]$.  It is an
  easy exercise to show that $U$ is Scott-open, and this concludes the
  argument.  \qed
\end{exa}

\begin{exa}
  \label{ex:dScott!=openball}
  There are cases where the $d$-Scott topology and the open ball
  topology are very different.  The most immediate examples are given
  by posets $X$, in which the $d_\leq$-Scott topology is the Scott
  topology, and the open ball topology is the Alexandroff topology,
  whose open sets are all the upwards-closed subsets.  \qed
\end{exa}

\begin{exa}
  \label{ex:dScott!=openball:Rl}
  For another example, consider the Sorgenfrey line $\Rl$.  We have
  seen that its open ball topology has a basis of half-open intervals
  $[x, x+r)$.  This is a bizarre topology: it is paracompact Hausdorff
  hence $T_4$ but the topological product of $\Rl$ with itself is not
  normal; it is first-countable but not countably-based; it is
  zero-dimensional, and not locally compact.  The $\dRl$-Scott
  topology is tamer: it has a basis of opens of the form
  $\Rl \cap \uuarrow (x, r)$, namely the open intervals $(x, x+r)$,
  hence it is just the usual topology on $\real$.  (We write
  $\uuarrow a$, in general, for $\{b \mid a \ll b\}$.)  \qed
\end{exa}

\section{Standard Quasi-Metric Spaces}
\label{sec:stand-quasi-metr}

The Scott topology, algebraicity and continuity, are mostly studied on
dcpos, not posets.  However, those notions do make sense on general
posets.  Similarly, research on quasi-metric spaces, their spaces of
formal balls, their generalized Scott topology, and so on, mostly
focused on Yoneda-complete spaces, but those should have a meaning
even in non-complete spaces.

However, non-complete spaces may exhibit a few pathologies that we
would like to exclude.  Here is how.  We shall use the notions and
results in subsequent sections.
\begin{defi}[Standard Quasi-Metric Space]
  \label{defn:d:std}
  \index{standard quasi-metric space}%
  A quasi-metric space $X, d$ is \emph{standard} if and only if, for
  every directed family of formal balls ${(x_i, r_i)}_{i \in I}$, for
  every $s \in \Rp$, ${(x_i, r_i)}_{i \in I}$ has a supremum in
  $\mathbf B (X, d)$ if and only if ${(x_i, r_i+s)}_{i \in I}$ has a
  supremum in $\mathbf B (X, d)$.
\end{defi}

Many quasi-metric spaces are standard, as we observe now.  Here and in
the sequel, a \emph{net} is a family ${(z_i)}_{i \in I, \sqsubseteq}$
of points $z_i$ indexed by a set $I$ equipped with a quasi-ordering
$\sqsubseteq$ that makes $I$ directed.  A limit of that net is any
point $z$ such that every open neighborhood $U$ of $z$ contains $z_i$
for $i$ large enough.
\begin{prop}
  \label{prop:d:std}
  Every metric space is standard.  Every Yoneda-complete quasi-metric
  space is standard.  Every poset is standard.
\end{prop}
\proof When $X, d$ is metric, we need to observe that: $(*)$ $(x, r)$
is the supremum of the directed family ${(x_i, r_i)}_{i \in I}$ if and
only if $x$ is the limit of the net ${(x_i)}_{i \in I, \sqsubseteq}$
(in $X$ with its open ball topology) and $r = \inf_{i \in I} r_i$
\cite[Theorem~5]{EH:comp:metric}.  The quasi-ordering $\sqsubseteq$ is
defined by $i \sqsubseteq j$ if and only if
$(x_i, r_i) \leq^{d^+} (x_j, r_j)$, and since that is directed,
${(x_i)}_{i \in I, \sqsubseteq}$ is a net.
A consequence of $(*)$ is that ${(x_i, r_i+s)}_{i \in I}$ has a
supremum $(x', r')$ in $\mathbf B (X, d)$ if and only if $x'$ is the
limit of ${(x_i)}_{i \in I, \sqsubseteq}$ and
$r' = \inf_{i \in I} r_i+s$.  It follows that the existence of a
supremum is equivalent for ${(x_i, r_i)}_{i \in I}$ and for
${(x_i, r_i+s)}_{i \in I}$, both being equivalent to the existence of
a limit of the net ${(x_i)}_{i \in I, \sqsubseteq}$.

When $X, d$ is Yoneda-complete, we use Lemma~7.4.25 and Lemma~7.4.26
of \cite{JGL-topology}, which together state the similar result that
$(x, r)$ is the supremum of the directed family
${(x_i, r_i)}_{i \in I}$ if and only if $x$ is the $d$-limit of
${(x_i)}_{i \in I, \sqsubseteq}$ and $r = \inf_{i \in I} r_i$.  The
notion of $d$-limit is irrelevant for our purposes here.  (We shall
need it later, and we shall define what it is there.)  The important
point is that, as above, the existence of a supremum is equivalent for
${(x_i, r_i)}_{i \in I}$ and for ${(x_i, r_i+s)}_{i \in I}$, both
being equivalent to the existence of a $d$-limit of the net
${(x_i)}_{i \in I, \sqsubseteq}$.

If $X$ is a poset, we have seen in Example~\ref{ex:<=} that
$(x, r) \mapsto (x, -r)$ defines an order-isomorphism from
$\mathbf B (X, d_\leq)$ onto $X \times (-\infty, 0]$.  Since suprema
in the latter are taken componentwise, and since all suprema exist in
$(-\infty, 0]$, the last claim is clear.  \qed

\begin{rem}
  \label{rem:d:std:no}
  Not every quasi-metric space is standard.  For a counterexample, let
  $X = [0, 1]$ with the quasi-metric defined by $d (x, y) = |x-y|$ if
  $x, y \neq 0$, $d (0, x) = a$ for every $x \neq 0$, and
  $d (x, 0) = 0$ for every $x$, where $a \geq 1$ is a fixed constant.
  (We need $a \geq 1$ to make sure the triangular inequality holds,
  notably $d (x, y) \leq d (x, 0) + d (0, y)$ for all $x, y$.)

  The directed family of formal balls $(1/2^m, 1/2^m)$, $m \in \nat$,
  has $(0, 0)$ has least upper bound, which one can show by verifying
  that $(0, 0)$ is in fast its sole upper bound.  Indeed, if $(x, r)$
  is an upper bound, and $x\neq 0$, then $|1/2^m - x| \leq 1/2^m - r$
  for every $m$, and as $m$ tends to $+\infty$, this forces $x=r=0$, a
  contradiction.  If $x=0$, then $d(1/2^m, x) = 0 \leq 1/2^m - r$,
  which forces $r=0$.

  For $s > 0$, the upper bounds of the directed family
  $(1/2^m, 1/2^m + s)$, $m \in \nat$, are those formal balls $(x, r)$
  such that $x+r \leq s$.  This proceeds in the same way as above.  If
  $(x,r)$ is an upper bound, and $x\neq 0$, then
  $|1/2^m - x| \leq 1/2^m+s-r$ for every $m$, and as $m$ tends to
  $+\infty$, $x \leq s-r$.  Conversely, if $x\leq s-r$, then
  $|1/2^m - x| \leq 1/2^m + x \leq 1/2^m+s-r$, so $(x, r)$ is an upper
  bound of the family.  If $x=0$, then $d(1/2^m, x)=0 \leq 1/2^m+s-r$
  for every $m$ implies $s\geq r$, which is again equivalent to
  $x +r \leq s$.  Conversely, if $s\geq r$ then
  $d (1/2^m, x) = 0 \leq 1/2^m +s-r$ for every $m$.

  For example, $(s/3, 2s/3)$ is such an upper bound.  We shall see in
  Proposition~\ref{prop:d:std:props}~(2) that, if $X, d$ were
  standard, then the least such upper bound would be $(0, s)$.  Now
  pick $s$ so that $0 < s < 3$: then $(0, s) \not\leq (s/3, 2s/3)$,
  since $d (0, s/3) = a \not\leq s-2s/3$.  This shows that $(0, s)$ is
  not least among all upper bounds of the family, hence that $X, d$
  cannot be standard.
\end{rem}

\begin{prop}
  \label{prop:d:std:props}
  In a standard quasi-metric space $X, d$, the following hold:
  \begin{enumerate}
  \item for every directed family of formal balls ${(x_i, r_i)}_{i \in
      I}$ with supremum $(x, r)$, $r = \inf_{i \in I} r_i$;
  \item for every directed family of formal balls
    ${(x_i, r_i)}_{i \in I}$ with supremum $(x, r)$, for every
    $s \in \real$ such that $s \geq -r$, the supremum of
    ${(x_i, r_i+s)}_{i \in I}$ exists and is equal to $(x, r+s)$;
    \index{radius map}%
  \item the \emph{radius map} $(x, r) \mapsto r$ is Scott-continuous
    from $\mathbf B (X, d)$ to $\Rp^{op}$ (the set of non-negative real
    numbers with the opposite ordering $\geq$);
  \item the map $\_ + s \colon (x, r) \mapsto (x, r+s)$ is
    Scott-continuous from $\mathbf B (X, d)$ to itself.
  \end{enumerate}
\end{prop}
\proof (1) Let $r_\infty = \inf_{i \in I} r_i$.  Observe that, for all
formal balls, $(y, s) \leq^{d^+} (z, t)$ implies $s \geq t$.  Since
$(x_i, r_i) \leq^{d^+} (x, r)$, $r_i \geq r$ for every $i \in I$, so
$r_\infty \geq r$.  By Definition~\ref{defn:d:std} with
$s = r_\infty - r$, the family ${(x_i, r_i - r_\infty)}_{i \in I}$
also has a supremum $(x', r')$, and by similar reasoning
$r_i - r_\infty \geq r'$ for every $i$.  This implies
$0 = \inf_{i \in I} r_i - r_\infty \geq r'$, so $r' = 0$.  In
particular, $(x', 0)$ is an upper bound of
${(x_i, r_i - r_\infty)}_{i \in I}$, namely,
$d (x_i, x') \leq r_i - r_\infty$ for every $i \in I$.  Equivalently,
$(x_i, r_i) \leq^{d^+} (x', r_\infty)$.  It follows that
$(x', r_\infty)$ is an upper bound of ${(x_i, r_i)}_{i \in I}$, hence
is above the least one, $(x, r)$.  In particular, $r \geq r_\infty$.
We have already proved the converse inequality, so $r = r_\infty$.

(2) By (1), $r = \inf_{i \in I} r_i$.  By Definition~\ref{defn:d:std},
${(x_i, r_i+s)}_{i \in I}$ has a supremum $(x', r')$, and by (1) again,
$r' = \inf_{i \in I} r_i + s = r + s$.  We now use the fact that
$(x_i, r_i) \leq^{d^+} (x, r)$ for every $i \in I$, i.e.,
$d (x_i, x) \leq r_i - r = (r_i+s) - (r+s)$, from which it follows
$(x_i, r_i+s) \leq^{d^+} (x, r+s)$ for every $i \in I$.  The formal
ball $(x, r+s)$ is an upper bound of ${(x_i, r_i+s)}_{i \in I}$, hence
is above its least upper bound $(x', r+s)$:
$d (x', x) \leq (r+s) - (r+s) = 0$.  Working in the converse
direction, the formal ball $(x', r) = (x', r+s + (-s))$ is an upper
bound of ${(x_i, r_i+s+(-s))}_{i \in I}= {(x_i, r_i)}_{i \in I}$,
hence is above $(x, r)$, so that $d (x, x') \leq r-r=0$.  Because
$X, d$ is quasi-metric, and $d (x, x') = d (x', x) = 0$, $x=x'$.

(3) The radius map is monotonic, namely $(x, r) \leq^{d^+} (y, s)$
implies $r \geq s$ (recall $\Rp^{op}$ has the opposite ordering
$\geq$), and what remains to be shown is (1).

(4) If $(x, r) \leq^{d^+} (x', r')$, then
$d (x, x') \leq r-r' = (r+s) - (r'+s)$, so
$(x, r+s) \leq^{d^+} (x', r')$.  This shows that $\_ + s$ is monotone.
Scott-continuity per se follows from (2).  \qed

The mapping $x \mapsto (x, 0)$ allows us to see $X$ (with the
$d$-Scott topology) as a topological subspace of $\mathbf B (X, d)$
(with the Scott topology).  If $(x, 0) \leq^{d^+} (y, s)$, then $0
\geq s$, so $s=0$.  It follows:
\begin{fact}
  \label{fact:X:upper}
  Let $X, d$ be a quasi-metric space.  Then $X$ embeds as an
  upwards-closed subset of $\mathbf B (X, d)$.
\end{fact}

When $X, d$ is standard, we can say more.
\begin{prop}
  \label{prop:X:Gdelta}
  Let $X, d$ be a standard quasi-metric space.  Then $X$ embeds as a
  $G_\delta$ subset of $\mathbf B (X, d)$.
\end{prop}
\proof Let $U_n = \{(x, r) \in \mathbf B (X, d) \mid r < 1/2^n\}$.
This is the inverse image of $[0, 1/2^n)$ by the radius map.  Since
the former is open in the Scott topology of $\Rp^{op}$, and using
Proposition~\ref{prop:d:std:props}~(3), $U_n$ is open.  Clearly, $X =
\bigcap_{n \in \nat} U_n$.  \qed

\section{Continuous Quasi-Metric Spaces}
\label{sec:cont-yoneda-compl}

In domain theory, there are dcpos, continuous dcpos, and algebraic
dcpos.  There are quasi-metric analogies of each notion, and we have
described them, except for continuous dcpos.

The definition of a \emph{continuous} Yoneda-complete quasi-metric
space stems from enriched category-theoretic considerations, and is
pretty complicated.  The first author claimed in
\cite[Definition~7.4.72]{JGL-topology} that a quasi-metric space
$X, d$ is continuous Yoneda-complete, in the sense of Kostanek and
Waszkiewicz, if and only $\mathbf B (X, d)$ is a continuous dcpo.
That happens to be true, as we shall see, but no proof is given of
that claim there.  We repair this omission, and also deal not only
with Yoneda-complete quasi-metric spaces, but with standard
quasi-metric spaces.

In order to define continuous Yoneda-complete quasi-metric spaces, in
principle, we need to first define the way-below $\mathcal Q$-relation
$\mathbf w \colon X \times X \to \creal$.  This is defined by an
enriched category-theoretic analogue of the definition of the usual
way-below relation on posets.  Fortunately, we will not need the
actual definition.  Kostanek and Waszkiewicz show that
\cite[Theorem~9.1]{KW:formal:ball}:
\begin{itemize}
\item If $X, d$ is continuous Yoneda-complete then $\mathbf B
  (X, d)$ is a continuous dcpo with way-below relation given by $(x,
  r) \ll (y, s)$ if and only if $r > \mathbf w (x, y) + s$.
\item If, for some map $v \colon X \times X \to \creal$,
  $\mathbf B (X, d)$ is a continuous dcpo and its way-below relation
  is characterized by $(x, r) \ll (y, s)$ if and only if
  $r > v (x, y) + s$, then $X, d$ is a continuous Yoneda-complete
  quasi-metric space.
\end{itemize}

We propose a simpler characterization, as a first step towards our
final simplification.  The key notion is a new twist on standardness.
\begin{defi}
  \label{defn:std}
  A binary relation $R$ on $\mathbf B (X, d)$ is \emph{standard} if and
  only if, for every $a \in \Rp$, $(x, r) \mathrel{R} (y, s)$ if and
  only if $(x, r+a) \mathrel{R} (y, s+a)$.
\end{defi}
This can be generalized to relations of any arity, including infinite
arities.  This way, this new notion of standardness also encompasses
Definition~\ref{defn:d:std}.  Note that $\leq^{d^+}$ is always standard.

We are interested in the cases where $\ll$ is standard.  Half of the
equivalence defining standardness is automatic:
\begin{lem}
  \label{lemma:ll:std:half}
  Let $X, d$ be a standard quasi-metric space.  For every $a \in \Rp$,
  if $(x, r+a) \ll (y, s+a)$ then $(x, r) \ll (y, s)$.
\end{lem}
\proof Let ${(z_i, t_i)}_{i \in I}$ be a directed family of formal
balls with a supremum $(z, t)$ above $(y, s)$.  By
Proposition~\ref{prop:d:std:props}~(2), ${(z_i, t_i+a)}_{i \in I}$ is
a directed family, and its supremum $(z, t+a)$ is above $(y, s+a)$.
Hence there is an index $i$ such that
$(y_i, s_i+a) \leq^{d^+} (z, t+a)$.  It follows that
$(y_i, s_i) \leq^{d^+} (z, t)$.  \qed


The converse implication, namely that $(x, r) \ll (y, s)$ implies
$(x, r+a) \ll (y, s+a)$, is wrong in general, as the following example
shows.

\begin{figure}
  \centering
  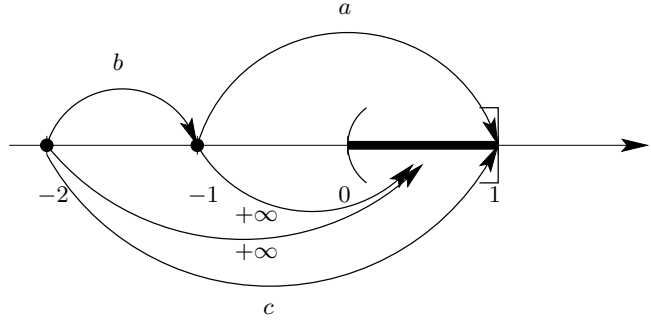
  \caption{A standard quasi-metric space with a non-standard way-below
    relation on formal balls}
  \label{fig:ng:nonstd}
\end{figure}
\begin{exa} 
  \label{ex:ng:nonstd}
  Let $X$ be the disjoint union of $(0, 1]$ with the Sorgenfrey
  quasi-metric $\dRl$, with two extra points which we shall name $-2$
  and $-1$.  We extend $\dRl$ to a quasi-metric $d$ on $X$ by letting
  $d (x, y) = +\infty$ if $x > y$, and when $x \leq y$ we let:
  $d (x, y) = y-x$ if $x, y \in (0, 1]$, $d (-1, y)=+\infty$ if
  $y \in (0, 1)$, $d (-1, 1)=a$, $d (-2, -1)=b$, $d (-2, y)=+\infty$
  if $y \in (0, 1)$, $d (-2, 1)=c$, where $a$, $b$ and $c$ are
  elements of $\Rp$ such that $c \leq a+b$ and $a, b > 0$.  See
  Figure~\ref{fig:ng:nonstd}, where only the distances from $x$ to $y$
  with $x < y$ are depicted.  To check that it is a quasi-metric, note
  that $d (x,z) \leq d (x,y) + d (y,z)$ is trivial whenever the
  right-hand side is equal to $+\infty$, or when any two points from
  $x$, $y$, $z$ are equal; so the only cases we have to check are
  those where $x < y < z$, and when $x=-2$, $y=-1$, $z=1$, this
  requires the inequality $c \leq a+b$.

  Given any directed family ${(x_i, r_i)}_{i \in I}$ of formal balls,
  observe that $(x_i, r_i) \leq^{d^+} (x_j, r_j)$ forces
  $x_i \leq x_j$, since otherwise
  $d (x_i, x_j) = +\infty \not\leq r_i - r_j$.  Then, either some
  $x_i$ is in $(0, 1]$, and we see that the cofinal family of those
  elements $(x_i, r_i)$ with $x_i \in (0, 1]$ has a supremum which is
  given as in $\Rl$, namely $(x, r)$ where $r = \inf_{i \in I} r_i$
  and $x = \sup_{i \in } x_i$; or ${(x_i, r_i)}_{i \in I}$ is included
  in $\{-2, -1\}$ and some $x_i$ equals $-1$, in which case there is a
  cofinal family where each $x_i$ is equal to $-1$, so the supremum
  is $(-1, \inf_{i \in I}, r_i)$; or every $x_i$ is equal to $-2$, and
  the supremum is $(-2, \inf_{i \in I} r_i)$.  In any case, the family
  has a supremum, so $\mathbf B (X, d)$ is a dcpo, in other words
  $X, d$ is Yoneda-complete.  In particular, it is standard.

  However, $\ll$ is not standard.  To this end, we fix some arbitrary
  real number $b' > b$, and we check that $(-2, b') \ll (-1, 0)$.
  Assume a monotone net ${(x_i, r_i)}_{i \in I, \sqsubseteq}$ whose
  supremum $(x, r)$ is above $(-1, 0)$.  Since $d (-1, x) \leq 0-r$,
  we must have $r=0$, and then a case analysis on $x$ shows that
  $x=-1$.  (The inequality $a > 0$ serves to show that $x=1$ is
  impossible.)  Since $r = \inf_{i \in I} r_i$, $r_i$ must be strictly
  less than $b > 0$ for $i$ large enough, and then the inequality
  $(x_i, r_i) \leq^{d^+} (x, r) = (-1, 0)$ implies $x_i=-1$.  Then
  $d (-2, -1) = b \leq b'-r_i$ for $i$ large enough, which shows that
  $(-2, b') \leq^{d^+} (-1, r_i) = (x_i, r_i)$ for $i$ large enough.

  It remains to check that $(-2, b'+a)$ is not way-below $(-1, a)$.
  Let $(x_n, r_n) = (1-1/2^{n+1}, 1/2^{n+1})$, a directed family in
  $X$ included in $(0, 1]$.  Its supremum is $(x, r) = (1, 0)$, which
  is above $(-1, a)$ since $d (-1, 1) = a$.  But no $(x_n, r_n)$ is
  above $(-2, b'+a)$, since that would mean that
  $d (-2, 1-1/2^{n+1}) = +\infty$ would be at most $b'+a-1/2^{n+1}$,
  which is impossible.  \qed
\end{exa}

We need the following easy lemma for the next theorem, and for later
results as well.

\begin{lem}
  \label{lemma:+eps}
  Let $X, d$ be a quasi-metric space.  For every
  $(y, s) \in \mathbf B (X, d)$, $(y, s)$ is the supremum in
  $\mathbf B (X, d)$ of the chain $(y, s+1/2^n)$, $n \in \nat$.
\end{lem}
\proof Clearly $(y, s+1/2^n) \leq^{d^+} (y, s)$ for every
$n \in \nat$.  For every upper bound $(z, t)$ of the elements
$(y, s+1/2^n)$, $n \in \nat$, $d (y, z) \leq s+1/2^n-t$ for every
$n \in \nat$, so $d (y, z) \leq s-t$, showing that
$(y, s) \leq^{d^+} (z, t)$.  \qed

\begin{prop}
  \label{prop:cont}
  A quasi-metric space $X, d$ is a continuous Yoneda-complete
  quasi-metric space if and only if $\mathbf B (X, d)$ is a continuous
  dcpo with a standard way-below relation.
\end{prop}
\proof If $X, d$ is continuous Yoneda-complete, then
$(x, r) \ll (y, s)$ if and only if $r < \mathbf w (x, y) + s$.  This
is clearly a standard relation.

Conversely, assume that $\mathbf B (X, d)$ is a continuous dcpo and
that $\ll$ is standard.  Let
$v (x, y) = \inf \{r-s \mid (x, r) \ll (y, s)\}$, the infimum being
equal to $+\infty$ if the right-hand set is empty.  We claim that
$(x, r) \ll (y, s)$ if and only if $r > v (x,y) + s$.

In one direction, assume $(x, r) \ll (y, s)$.  By interpolation, find
a formal ball $(z, t)$ such that $(x, r) \ll (z, t) \ll (y, s)$.
Using Lemma~\ref{lemma:+eps}, $(z, t) \leq^{d^+} (y, s+1/2^n)$ for
some $n \in \nat$, hence $(x, r) \ll (y, s+1/2^n)$.  By definition of
$v$, $v (x, y) \leq r - s - 1/2^n$, so that $v (x, y) < r-s$.

Conversely, if $r > v (x,y) +s$, then by definition of $v$, there are
numbers $r', s' \in \real^+$ such that $(x, r') \ll (y, s')$ and
$r > r'-s'+s$.  By Lemma~\ref{lemma:ll:std:half},
$(x, r'-s') \ll (y, 0)$, and since $\ll$ is standard,
$(x, r'-s'+s) \ll (y, s)$.  Using $r > r'-s'+s$, we obtain that
$(x, r) \leq^{d^+} (x, r'-s'+s) \ll (y, s)$.  \qed

We now refine this by showing that the continuity of
$\mathbf B (X, d)$ (for standard $X, d$) is enough to ensure that
pathologies such as Example~\ref{ex:ng:nonstd} do not actually happen.
\begin{prop}
  \label{prop:ng}
  Let $X, d$ be a standard quasi-metric space.  If $\mathbf B (X, d)$
  is a continuous poset, then it has a standard way-below relation.
\end{prop}
\proof Let $(x, r) \ll (y, s)$, and fix $a \in \Rp$.  The family
$\ddarrow (y, s+a)$ is directed, has $(y, s+a)$ as supremum, and
consists of elements whose radius is at least $s+a \geq a$.  Write
that family ${(z_i, t_i+a)}_{i \in I}$.  Note that
$(z_i, t_i+a) \leq^{d^+} (z_j, t_j+a)$ if and only if
$(z_i, t_i) \leq^{d^+} (z_j, t_j)$, because both are equivalent to
$d (z_i, z_j) \leq t_i - t_j$.  Therefore ${(z_i, t_i)}_{i \in I}$ is
also a directed family.  By Proposition~\ref{prop:d:std:props}~(2), it
admits $(y, s)$ as supremum.  Since $(x, r) \ll (y, s)$,
$(x, r) \leq^{d^+} (z_i, t_i)$ for some $i \in I$.  It is easy to see
that $(x, r+a) \leq^{d^+} (z_i, t_i+a)$.  Since
$(z_i, t_i+a) \in \ddarrow (y, s+a)$, $(x, r+a) \ll (y, s+a)$.  \qed
Together with Proposition~\ref{prop:cont}, we therefore obtain:
\begin{thm}
  \label{thm:cont}
  A quasi-metric space $X, d$ is a continuous Yoneda-complete
  quasi-metric space if and only if $\mathbf B (X, d)$ is a continuous
  dcpo.  \qed
\end{thm}

\begin{exa}
  \label{ex:cont:R}
  Following up on Example~\ref{ex:Yoneda-complete:R}, $\creal, \dreal$
  is continuous Yoneda-complete.
  It is even algebraic Yoneda-complete, as we shall see in
  Example~\ref{ex:alg:R}.  \qed
\end{exa}

\begin{exa}
  \label{ex:cont:Rl}
  $\Rl$ is continuous Yoneda-complete
  (see Example~\ref{ex:Yoneda-complete:Rl}).
  \qed
\end{exa}

Much as continuous dcpos can be generalized to continuous posets, this
allows us to define continuous quasi-metric spaces without
Yoneda continuity.
\begin{defi}
  \label{defn:cont}
  A standard quasi-metric space $X, d$ is \emph{continuous} if and
  only if $\mathbf B (X, d)$ is a continuous poset.
\end{defi}

\begin{exa}
  \label{ex:cont:metric}
  Edalat and Heckmann noticed that, when $X, d$ is a metric space,
  $\mathbf B (X, d)$ is always a continuous poset, with
  $(x, r) \ll (y, s)$ if and only if $d (x, y) < r-s$
  \cite{EH:comp:metric}.  It follows that all metric spaces are
  continuous.  
  \qed
\end{exa}

\begin{exa}
  \label{ex:cont:<=}
  Let us return to the case of posets (Example~\ref{ex:<=}).  For a
  poset $X$, recall that $\mathbf B (X, d_\leq)$ is order-isomorphic
  to $X \times (-\infty, 0]$.  The latter is continuous if and only if
  $X$ is continuous.
  Hence a poset is continuous qua quasi-metric space if and only if it
  is continuous in the usual sense.  \qed
\end{exa}

\section{Continuous Yoneda-complete Spaces are Sober, Choquet-complete}
\label{sec:choquet-completeness}

We now observe that continuous Yoneda-complete quasi-metric spaces have a
number of desirable properties: they are sober, and they are
Choquet-complete, in particular they are Baire.

\begin{prop}
  \label{prop:cont:sober}
  Every continuous Yoneda-complete quasi-metric space $X, d$ is sober
  in its $d$-Scott topology.
\end{prop}
\proof Since $\mathbf B (X, d)$ is a continuous dcpo, it is sober in
its Scott topology, see e.g.\ \cite[Proposition~7.2.27]{AJ:domains},
\cite[Corollary~II-1.12]{GHKLMS:contlatt}, or
\cite[Proposition~8.2.12~$(b)$]{JGL-topology}.  Consider the diagram:
\[
\xymatrix{
  \mathbf B (X, d) \ar@<1ex>[r]^{rad} \ar@<-1ex>[r]_0 & \Rp^{op}
}
\]
where $rad$ is the radius map and $0$ is the constant $0$ map.  Both
$rad$ and $0$ are continuous, using Proposition~\ref{prop:d:std:props}~(3).
We claim that the map $\eta \colon x \in X \mapsto (x, 0)$ is an equalizer of that
diagram in the category of topological spaces.  Consider any
continuous map $f \colon Y \to \mathbf B (X, d)$ such that $rad \circ
f = 0 \circ f$.  For every $y \in Y$, $f (y)$ is an element of the
form $(x, 0)$, and we define $f' (y)$ as $x$.  Clearly, $f = \eta
\circ f'$.  Moreover, $f'$ is continuous: for every open subset $U$ of
$X$, by definition of the $d$-Scott topology there is a Scott-open
subset $V$ of $\mathbf B (X, d)$ such that $U = \eta^{-1} (V)$, and
then ${f'}^{-1} (U) = f^{-1} (V)$ is open.  The result follows because
spaces obtained as equalizers of continuous maps from a sober space to
a topological space are sober; see, e.g.,
\cite[Lemma~8.4.12]{JGL-topology}.  \qed

Let us turn to Choquet completeness.

Given a topological space $X$, the \emph{strong Choquet game} on $X$
is played as follows.  There are two players, $\alpha$ and $\beta$,
who alternate turns.  Player $\beta$ starts and chooses a non-empty
open subset $V_0$ of $X$, and a point $x_0 \in V_0$.  Then $\alpha$
plays an open subset $U_0$ of $V_0$ containing $x_0$.  Player $\beta$
finds a non-empty open subset $V_1$ of $U_0$, and picks a point
$x_1 \in V_1$, then $\alpha$ produces an open subset $U_1$ of $V_1$
containing $x_1$, and so on.  Clearly,
$\bigcap_{n \in \nat} U_n = \bigcap_{n \in \nat} V_n$, and we say that
$\alpha$ wins the game if and only if that set is non-empty.  A
strategy for $\alpha$ is a map from \emph{histories}
$x_0, V_0, U_0, x_1, V_1, U_1, \cdots, x_n, V_n$ (namely, $n\in \nat$,
all $U_i$ and $V_i$ are open,
$V_0 \supseteq U_0 \supseteq \cdots \supseteq U_{n-1} \supseteq V_n$,
$x_0 \in U_0$, \ldots, $x_{n-1} \in U_{n-1}$, $x_n \in V_n$) to opens
$U_n$ such that $x_n \in U_n \subseteq V_n$.  $X$ is
\emph{Choquet-complete} if and only if $\alpha$ has a winning
strategy.  See \cite[Section~7.6]{JGL-topology}.

Every continuous dcpo is Choquet-complete in its Scott topology, an
observation due to K. Martin \cite{Martin:nonclassical}.   Player
$\alpha$'s winning strategy can even be chosen to be \emph{stationary}
\cite[Lemma~7.6.3]{JGL-topology}, i.e., so that $U_n$ depends only on
$\beta$'s last move $x_n, V_n$; and \emph{convergent}
\cite{DM:choquet}, i.e., so that ${(U_n)}_{n \in \nat}$, or
equivalently ${(V_n)}_{n \in \nat}$, is a neighborhood base of some
element $y$.  (We call a space \emph{convergence Choquet-complete}
if and only if $\alpha$ has a convergent winning strategy.)
The argument is simple: given $\beta$'s last move $x_n,
V_n$, $\alpha$ picks an element $y_n \ll x_n$ such that $y_n \in V_n$,
and plays $U_n = \uuarrow y_n$; then $y = \sup_{n \in \nat} y_n$.

Every $G_\delta$ subset of a Choquet-complete space is Choquet
complete.  This is mentioned as Theorem~2.30~$(iv)$ in
\cite{Martin:nonclassical}, and a proof can be found in
\cite[Proposition~2.1~$(iii)$]{HKL:glimmeffros}.
The same proof shows:
\begin{lem}
  \label{lemma:Gdelta:Choquet}
  Every $G_\delta$ subset of a convergence Choquet-complete space is a
  convergence Choquet-complete subspace.
\end{lem}
\proof Let $G$ be a $G_\delta$ subset of a convergence
Choquet-complete space $X$, and write $G$ as the intersection of a
decreasing sequence of opens
$G_0 \supseteq G_1 \supseteq \cdots \supseteq G_n \supseteq \cdots$.
For convenience, we assume $G_0 = X$.  Every open subset $V$ of $G$
can be written as $V' \cap G$ for some open subset $V'$ of $X$.
Picking the largest such open for $V'$, we can ensure that the
assignment $V \mapsto V'$ is monotonic.

Assume a convergent winning strategy $\sigma$ for $\alpha$ on $X$.  We
obtain a winning strategy for $\alpha$ on $G$ as follows: on the
history $x_0, V_0, U_0, x_1, V_1, U_1, \cdots, x_n, V_n$ (inside $G$),
$\alpha$ computes $U_n = G \cap W_n$, where
$W_n = G_n \cap \sigma (x_0, V'_0, U'_0 \cap G_0, x_1, V'_1, U'_1 \cap
G_1, \cdots, x_n, V'_n)$,
therefore simulating a play on $X$ with strategy $\sigma$.  (Note that
each $x_i$, $1\leq i \leq n-1$, is indeed in $U'_{i-1} \cap G_{i-1}$,
since $x_i \in U_i = U'_i \cap G$, and $x_n \in V'_n$, $x_n \in G_n$.
Note also that $W_n = U'_n \cap G_n$.)  By assumption,
${(V'_n)}_{n \in \nat}$ is a neighborhood base of some point
$y \in X$.  By construction,
$y \in \bigcap_{n \in \nat} W_n \subseteq \bigcap_{n \in \nat} G_n =
G$.
Every open neighborhood $V$ of $y$ in $G$ is such that
$V'_n \subseteq V'$ for some $n \in \nat$, so
$V_n = V'_n \cap G \subseteq V' \cap G = V$, showing that
${(V_n)}_{n \in \nat}$ is a neighborhood base of $y$ in $G$.  \qed


Recalling Proposition~\ref{prop:X:Gdelta} and the fact that every
Yoneda-complete space is standard, we obtain the following.
\begin{thm}
  \label{thm:cont:Choquet}
  Every continuous Yoneda-complete quasi-metric space is convergence
  Choquet-complete in its $d$-Scott topology.
\end{thm}

\begin{rem}
  Theorem~\ref{thm:cont:Choquet} generalizes the fact that every
  Smyth-complete space is convergence Choquet-complete in its open
  ball topology \cite[Exercise~7.6.5]{JGL-topology}.  Indeed, every
  Smyth-complete space is continuous Yoneda-complete, and its open
  ball topology coincides with its $d$-Scott topology.
\end{rem}

Recall that every Choquet-complete space is Baire (see, e.g.,
\cite[Theorem~7.6.8]{JGL-topology}), namely, the intersection of
countably many dense open subsets is dense.
\begin{cor}
  \label{cor:cont:Baire}
  Every continuous Yoneda-complete quasi-metric space is Baire in its
  $d$-Scott topology.
\end{cor}



\section{Algebraic Quasi-Metric Spaces}
\label{sec:algebr-yoneda}

The original definition of \emph{algebraic} Yoneda-complete spaces is
pretty complicated, and the point we would like to make here is that
there is a simpler one, which extends naturally to non-complete spaces
as well.

This is the point where we have to recapitulate the standard
definitions.  Fix a quasi-metric space $X, d$.  A net
${(x_i)}_{i \in I, \sqsubseteq}$ is \emph{Cauchy} if and only if for
every $\epsilon > 0$, there is an $i_0 \in I$ such that for all
$i, j \in I$ with $i_0 \sqsubseteq i \sqsubseteq j$,
$d (x_i, x_j) < \epsilon$.  A point $x \in X$ is the \emph{$d$-limit}
of the Cauchy net ${(x_i)}_{i \in I, \sqsubseteq}$ if and only if, for
every $y \in X$, $d (x, y)=\limsup_{i \in I, \sqsubseteq} d (x_i, y)$.
A Yoneda-complete space is a quasi-metric space where every Cauchy net
has a $d$-limit.

The relation to formal balls is as follows.  Observe that, if
${(x_i, r_i)}_{i \in I, \sqsubseteq}$ is any monotone net of formal
balls such that $\inf_{i \in I} r_i = 0$, then ${(x_i)}_{i \in I}$ is
Cauchy \cite[Lemma~7.2.7]{JGL-topology}.  Such a monotone net is
called a \emph{Cauchy-weighted} net, as the numbers $r_i$ act as
weights that witness the fact that ${(x_i)}_{i \in I}$ is Cauchy.  A
net ${(x_i)}_{i \in I, \sqsubseteq}$ is \emph{Cauchy-weightable} if
and only if one can find weights $r_i$ that make
${(x_i, r_i)}_{i \in I, \sqsubseteq}$ Cauchy-weighted.  So every
Cauchy-weightable net is Cauchy.  The converse fails
\cite[Exercise~7.2.12]{JGL-topology}, but every Cauchy net has a
Cauchy-weightable subnet \cite[Lemma~7.2.8]{JGL-topology}, and they
behave similarly as far as $d$-limits are concerned: if a Cauchy net
has a $d$-limit, then all its subnets are Cauchy and have the same
$d$-limit \cite[Exercise~7.4.7]{JGL-topology}, and conversely, if a
Cauchy subnet of a Cauchy net has a $d$-limit, then this is also a
$d$-limit of the Cauchy net \cite[Lemma~7.4.6]{JGL-topology}.  In
particular, an equivalent definition of Yoneda-completeness is: every
Cauchy-\emph{weightable} net has a $d$-limit.  This trick finds its
roots in \cite[Section~2.2]{EH:comp:metric}.

\begin{rem}
  \label{rem:dlim:sup}
  While $d$-limits of Cauchy nets are defined through a limit
  superior, $d$-limits of Cauchy-weightable nets can be characterized
  by a simpler formula: given a Cauchy-weighted net
  ${(x_i, r_i)}_{i \in I, \sqsubseteq}$, $x$ is the $d$-limit of
  ${(x_i)}_{i \in I, \sqsubseteq}$ if and only if, for every
  $y \in X$, $d (x, y) = \sup_{i \in I} (d (x_i, y) - r_i)$
  \cite[Lemma~7.4.9]{JGL-topology}.  Moreover, the latter is a
  directed supremum.
\end{rem}

\begin{exa}
  \label{ex:Rp:dlim}
  In $\creal, \dreal$, every net ${(x_i)}_{i \in I, \sqsubseteq}$
  (even not Cauchy) has a $\dreal$-limit, which is its limit superior
  \cite[Exercise~7.1.16]{JGL-topology}.  Given a Cauchy-weighted net
  ${(x_i, r_i)}_{i \in I, \sqsubseteq}$, the $\dreal$-limit of
  ${(x_i)}_{i \in I, \sqsubseteq}$ can be expressed as the simpler,
  directed supremum $\sup_{i \in I} (x_i - r_i)$.  More generally, the
  supremum $(x, r)$ of any directed family ${(x_i, r_i)}_{i \in I}$ of
  formal balls in $\mathbf B (\creal, \dreal)$ is given by
  $r = \inf_{i \in I} r_i$, and $x$ is the directed supremum
  $\sup_{i \in I} (x_i + r - r_i)$.  \qed
\end{exa}

\begin{exa}
  \label{ex:poset:dlim}
  Look at the case of posets, seen as quasi-metric spaces.  A net
  ${(x_i)}_{i \in I, \sqsubseteq}$ is Cauchy if and only if it is an
  eventually monotone net, that is, for $i, j$ large enough,
  $i \sqsubseteq j$ implies $x_i \leq x_j$.  It is Cauchy-weightable
  if and only if it is a monotone net.  The notion of $d_\leq$-limit
  of Cauchy-weightable nets coincides with the notion of directed
  supremum.  \qed
\end{exa}

A point $x$ of a quasi-metric space $X, d$ is called \emph{$d$-finite}
if and only if, for every directed family ${(y_i, s_i)}_{i \in I}$ of
open balls with a supremum of the form $(y, 0)$, $d (x, y)$ is the
infimum of the filtered family ${(d (x, y_i) + s_i)}_{i \in I}$ of
elements of $\creal$ \cite[Lemma~7.4.56]{JGL-topology}.  This is not
the standard definition, which involves limits inferiors and Cauchy
nets, but it is closer to our needs.

A quasi-metric space $X, d$ is \emph{algebraic} if and only if every
point is a $d$-limit of some Cauchy net of $d$-finite points.  By the
same argument as above, it is equivalent to require that every point
be a $d$-limit of some Cauchy-\emph{weightable} net of $d$-finite
points.

We would like to offer a simpler view of those notions---at least on
standard quasi-metric spaces---based on formal balls.  The starting
point is the following.  We have two topologies on $\mathbf B (X, d)$,
the open ball topology of $d^+$, and the Scott topology of
$\leq^{d^+}$.  They have the same specialization ordering,
$\leq^{d^+}$, and the former is finer than the latter, as we see now.

\begin{lem}
  \label{lemma:Scott=>openball}
  Let $X, d$ be a quasi-metric space.  The open ball topology on
  $\mathbf B (X, d), d^+$ is finer than the Scott topology.
\end{lem}
\proof Let $U$ be a Scott-open subset of $\mathbf B (X, d)$, and
$(y, s)$ be a point in $U$.  By Lemma~\ref{lemma:+eps}, $(y, s+1/2^n)$
is in $U$ for some $n \in \nat$.  We check easily that $(y, s)$ is in the
open ball $B^{d^+}_{(y, s), <1/2^n}$
and that the latter open ball is included in $\upc (y, s+1/2^n)$,
which is included in $U$.
Hence we have found an open neighborhood of $(y, s)$ for the open ball
topology that is included in $U$.  It follows that $U$ is open in the
open ball topology.  \qed

Conversely, one may wonder when the open ball topology coincides with
the Scott topology on $\mathbf B (X, d)$.  We shall solve this
question below.

In the meantime, observe that we may separate the open balls in
$\mathbf B (X, d), d^+$ into the nice open balls, those which are
Scott-open, and those that are not nice.  We may then consider a point
$x \in X$ as nice if and only if all the open balls centered at
$(x, 0)$ are Scott-open.  Let us give a name to that notion.
\begin{defi}[Center point]
  \label{defn:d:center}
  \index{center point}%
  In a quasi-metric space $X, d$, an element $x \in X$ is a
  \emph{center point} if and only if, for every $\epsilon > 0$, the
  open ball $B^{d^+}_{(x, 0), <\epsilon}$ is Scott-open in
  $\mathbf B (X, d)$.
\end{defi}
The intersection of $B^{d^+}_{(x, 0), <\epsilon}$ with $X$ is just
$B^d_{x, <\epsilon}$, making the following immediate.
\begin{rem}
  \label{rem:d:center}
  For every center point $x$ of a hemi-metric space $X, d$, every open
  ball $B^d_{x, <\epsilon}$ centered at $x$ is open in the $d$-Scott
  topology.
\end{rem}

The following shows that we can replace the complex definition of
``$d$-finite'' by the more synthetic notion of center point, in
standard spaces.  The proof also shows that every center point is
$d$-finite, even without assuming standardness.
\begin{lem}
  \label{lemma:d:finite}
  Let $X, d$ be a standard quasi-metric space.  The following are
  equivalent, for every $x \in X$:
  \begin{enumerate}
  \item $x$ is $d$-finite;
  \item for all $r \in \Rp$ and $\epsilon > 0$, the open ball
    $B^{d^+}_{(x, r), <\epsilon}$ is Scott-open in $\mathbf B (X, d)$;
  \item $x$ is a center point.
  \end{enumerate}
\end{lem}
A related result was given in \cite[Theorem~8.3]{KW:formal:ball}, who
instead show that all points are $d$-finite if and only if all points
are center points.

\proof $(3) \limp (1)$.  Let ${(y_i, s_i)}_{i \in I}$ be a directed family
of open balls with supremum $(y, 0)$.  Since
$(y_i, s_i) \leq^{d^+} (y, 0)$ for every $i \in I$,
$d (y_i, y) \leq s_i$, so $d (x, y) \leq d (x, y_i) + s_i$, by the
triangular inequality.  Assume the inequality
$d (x, y) \leq \inf_{i \in I} (d (x, y_i) + s_i)$ were strict:
$d (x, y) < \epsilon \leq \inf_{i \in I} (d (x, y_i) + s_i)$ for some
$\epsilon > 0$.  Then $(y, 0)$ would be in
$B^{d^+}_{(x, 0), <\epsilon}$, so some $(y_i, s_i)$ would be in
$B^{d^+}_{(x, 0), <\epsilon}$, which is Scott-open by assumption.
That is, $d (x, y_i) +s_i < \epsilon$.  However, that would contradict
the inequality $\epsilon \leq \inf_{i \in I} (d (x, y_i) + s_i)$.
Therefore $d (x, y) = \inf_{i \in I} (d (x, y_i) + s_i)$, showing that
$x$ is $d$-finite.

$(1) \limp (2)$, assuming $X, d$ standard.  Let $x$ be a $d$-finite point
of $X$.  The open ball $B^{d^+}_{(x, r), <\epsilon}$ is the set of
formal balls $(y, s)$ such that $d^+ ((x, r), (y, s)) < \epsilon$,
i.e., such that $d (x, y) - r + s < \epsilon$.  This is upwards-closed
with respect to $\leq^{d^+}$, by the triangular inequality.  Let
$(y, s)$ be the supremum of some directed family
${(y_i, s_i)}_{i \in I}$ of formal balls, and assume that $(y, s)$ is
in $B^{d^+}_{(x, r), <\epsilon}$.  By
Proposition~\ref{prop:d:std:props}~(2),
$(y, 0) = \sup_{i \in I} (y_i, s_i - s)$, and the fact that $x$ is
$d$-finite then implies that
$d (x, y) = \inf_{i \in I} (d (x, y_i) + s_i -s)$.  Using the
inequality $d (x, y) -r+s < \epsilon$, there must exist an $i \in I$
such that $d (x, y_i) + s_i -r < \epsilon$, i.e., such that
$(y_i, s_i)$ is in $B^{d^+}_{(x, r), <\epsilon}$.  In other words,
$B^{d^+}_{(x, r), <\epsilon}$ is Scott-open.

$(2) \limp (3)$ is trivial: take $r=0$.  \qed

Further characterizations follow, assuming $X, d$ continuous.
\begin{lem}
  \label{lemma:center}
  Let $X, d$ be a continuous quasi-metric space.  The following are
  equivalent, for every $x \in X$:
  \begin{enumerate}
  \item $x$ is a center point;
  \item $v (x, \_) = d (x, \_)$, where $v (x, y)$ is defined as $\inf
    \{r-s \mid (x, r) \ll (y, s)\}$;
  \item for every $\epsilon > 0$, $B^{d^+}_{(x, 0), <\epsilon} =
    \uuarrow (x, \epsilon)$.
  \end{enumerate}
\end{lem}
(The map $v$ was already used in the proof of Proposition~\ref{prop:cont}.)

\proof $(1) \limp 2$.  If $(x, r) \ll (y, s)$, then
$(x, r) \leq^{d^+} (y, s)$, so $d (x, y) \leq r-s$.  Taking infima,
$d (x, y) \leq v (x, y)$.  If the inequality were strict, there would
be a positive real $\epsilon$ such that
$d (x, y) < \epsilon \leq v (x, y)$.  In particular,
$(y, 0) \in B^{d^+}_{(x, 0), <\epsilon}$.  Notice that
$B^{d^+}_{(x, 0), <\epsilon} \subseteq \upc (x, \epsilon)$: every
element $(z, t)$ of $B^{d^+}_{(x, 0), <\epsilon}$ is such that
$d (x, z) + t < \epsilon$, so $(x, \epsilon) \leq^{d^+} (z, t)$.
Since $x$ is a center point, $B^{d^+}_{(x, 0), <\epsilon}$ is included
in the Scott interior of $\upc (x, \epsilon)$, and that is
$\uuarrow (x, \epsilon)$: that the Scott interior of $\upc a$ is
$\uuarrow a$ is true in every continuous poset
\cite[Proposition~5.1.35]{JGL-topology}.  We have therefore shown that
$(y, 0) \in \uuarrow (x, \epsilon)$, that is,
$(x, \epsilon) \ll (y, 0)$.  It follows that $v (x, y) \leq \epsilon$,
a contradiction.

$(2) \limp (3)$.  For every $(y, s) \in B^{d^+}_{(x, 0), <\epsilon}$,
$d (x, y) + s < \epsilon$, and using $(2)$, $v (x, y) + s < \epsilon$.
That implies $\epsilon > s$, and the existence of $r', s' \in \Rp$
such that $(x, r') \ll (y, s')$ and $r'-s' < \epsilon-s$.  Since $\ll$
is standard by Proposition~\ref{prop:ng},
$(x, \epsilon) \ll (y, s'+\epsilon-r')$.  Since $r'-s' < \epsilon-s$,
$s'+\epsilon-r' > s$ and in particular the formal ball
$(y, s'+\epsilon-r')$ makes sense; further,
$(y, s'+\epsilon-r') \leq^{d^+} (y, s)$, so
$(x, \epsilon) \ll (y, s)$.  It follows that
$B^{d^+}_{(x, 0), <\epsilon} \subseteq \uuarrow (x, \epsilon)$.  The
converse direction is simpler.  For every
$(y, s) \in \uuarrow (x, \epsilon)$, there is an $n \in \nat$ such
that $(x, \epsilon) \leq^{d^+} (y, s+1/2^n)$ by
Lemma~\ref{lemma:+eps}.  That implies
$d (x, y) \leq \epsilon - s - 1/2^n < \epsilon-s$, whence
$(y, s) \in B^{d^+}_{(x, 0), <\epsilon}$.

$(3) \limp (1)$ is obvious.  \qed

\begin{exa}
  \label{ex:alg:R:center}
  The $\dreal$-finite points of $\creal$ are known to be those that
  are different from $+\infty$.  Equivalently, let us check that they
  are the center points.  We already know from Example~\ref{ex:cont:R}
  that $v (x, \_) = \dreal (x, \_)$ if and only if $x \neq +\infty$.
  It is now enough to apply Lemma~\ref{lemma:center}.  \qed
\end{exa}

\begin{exa}
  \label{ex:alg:Rl:center}
  \emph{No} point of $\Rl$ is $\dRl$-finite
  \cite[Exercise~7.4.73]{JGL-topology}.  We check this here.
  Recalling Example~\ref{ex:cont:Rl}, $\Rl$ is continuous
  Yoneda-complete, and we see that $v (x, \_)$ is never equal to the
  map $\dRl (x, \_)$, since $v (x, x) = +\infty$ but
  $\dRl (x, x) = 0$.  \qed
\end{exa}

To capture algebraicity, we define the following notion.
\begin{defi}[Enough center points]
  \label{defn:d:enoughcenter}
  A quasi-metric space $X, d$ \emph{has enough center points} if and
  only if the open balls $B^{d^+}_{(x, 0), <\epsilon}$, with $x$ a
  center point and $\epsilon > 0$, form a base of the Scott topology
  on $\mathbf B (X, d)$.
\end{defi}
In concrete terms, $X, d$ has enough center points if and only if, for
every Scott-open subset $\mathcal U$ of $\mathbf B (X, d)$, for every
$(y, s) \in \mathcal U$, there is a center point $x$ and an
$\epsilon > 0$ such that
$(y, s) \in B^{d^+}_{(x, 0), <\epsilon} \subseteq \mathcal U$.  This
is what we shall use in proofs.

\begin{exa}
  \label{ex:alg:R}
  $\creal, \dreal$ has enough center points.  Indeed, recall that
  $\mathbf B (\creal, \dreal)$ is isomorphic to the closed set
  $C = \{(a, b) \in (\real \cup \{+\infty\}) \times (-\infty, 0] \mid
  a-b \geq 0\}$
  of Example~\ref{ex:Yoneda-complete:R}.  The latter is a continuous
  dcpo with $(a, b) \ll (a', b')$ if and only if $a < a'$ and
  $b < b'$, and therefore the points $(a, b) \in C$ with
  $a \neq +\infty$ form a basis.  Using the isomorphism, a basis of
  $\mathbf B (\creal, \dreal)$ consists of all the formal balls
  $(x, r)$ with $x \neq +\infty$, hence a base of the Scott topology
  is given by the sets $\uuarrow (x, r)$, $x \neq +\infty$.  Using
  Lemma~\ref{lemma:center}, the latter are the open balls
  $B^{d^+}_{(x, 0), <r}$, and $x$ is a center point, as we have seen
  in Example~\ref{ex:alg:R:center}.  \qed
\end{exa}

\begin{exa}
  \label{ex:alg:Rl}
  Since no point of $\Rl$ is finite, $\Rl$ is far from having enough
  center points.  We shall see that having enough center points is
  equivalent to being algebraic, and thus we retrieve the fact that
  $\Rl$ is not algebraic (\cite{KW:formal:ball}, see also
  \cite[Exercise~7.4.73]{JGL-topology}).  \qed
\end{exa}

\begin{prop}
  \label{prop:d:enoughcenter}
  Let $X, d$ be a quasi-metric space with enough center points.  The
  open balls $B^d_{x, <\epsilon}$ with $x$ a center point and
  $\epsilon > 0$ form a base of the $d$-Scott topology on $X$.  If
  $X, d$ is standard, then we can additional require that all the
  involved radii $\epsilon$ are less than some prescribed upper bound
  $b > 0$.
\end{prop}
\proof Given $x \in X$ and an open neighborhood $U$ of $x$, let
$\mathcal U$ be some Scott-open subset of $\mathbf B (X, d)$ whose
intersection with $X$ equals $U$.  By definition, there is an open
ball $B^{d^+}_{(a, 0), <\epsilon}$ that contains $(x, 0)$ and is
included in $\mathcal U$, where $a$ is a center point and
$\epsilon > 0$.  Taking intersections with $X$, $B^d_{a, <\epsilon}$
is an open neighborhood of $x$ that is included in $U$.

If $X, d$ is standard, then for $b \in \Rp$, $b > 0$, the set
$V_b = \{(y, s) \mid s < b\}$ is open, as the inverse image of
$[0, b)$ by the radius map, which is continuous by
Proposition~\ref{prop:d:std:props}~(3).  Then $\mathcal U \cap V_b$ is
also Scott-open and contains $(x, 0)$, so there is an open ball
$B^{d^+}_{(a, 0), <\epsilon}$, where $a$ is a center point, that
contains $(x, 0)$ and is included in $\mathcal U \cap V_b$.  As above,
we conclude that $B^d_{a, <\epsilon}$ is an open neighborhood of $x$
that is included in $U$.  Additionally, since
$(a, s) \in B^{d^+}_{(a, 0), <\epsilon} \subseteq V_b$ for every
$s < \epsilon$, $\epsilon$ is less than or equal to $b$.  \qed

We have already announced the following result, as part of the
Kostanek-Waszkiewicz theorem, assuming $X, d$ Yoneda-complete.  We now
show that this holds in all standard quasi-metric spaces, even not
Yoneda-complete.
\begin{lem}
  \label{lemma:d:limit:sup}
  Let $X, d$ be a standard quasi-metric space.  For every
  Cauchy-weighted net ${(x_i, r_i)}_{i \in I, \sqsubseteq}$, a point
  $x \in X$ is a $d$-limit of ${(x_i)}_{i \in I, \sqsubseteq}$ if and
  only if $(x, 0)$ is the supremum of the directed family
  ${(x_i, r_i)}_{i \in I}$ in $\mathbf B (X, d)$.
\end{lem}
\proof Assume $x$ is a $d$-limit of ${(x_i)}_{i \in I, \sqsubseteq}$.
For every formal ball $(y, s)$, $(y, s)$ is an upper bound of
${(x_i, r_i)}_{i \in I}$ if and only if $d (x_i, y) \leq r_i-s$ for
every $i \in I$, if and only if $s=0$ (using $\inf_{i \in I} r_i=0$)
and $d (x_i, y) - r_i \leq 0$ for every $i \in I$.  Since, by
Remark~\ref{rem:dlim:sup},
$d (x, y)=\sup_{i \in I} (d (x_i, y) - r_i)$, $(y, s)$ is an upper
bound of ${(x_i, r_i)}_{i \in I}$ if and only if $s=0$ and
$d (x, y) = 0$, and that is equivalent to $d (x, y) \leq 0 - s$, i.e.,
to $(x, 0) \leq^{d^+} (y, s)$.  The least such formal ball $(y, s)$ is
then $(x, 0)$.

Conversely, assume that $(x, 0)$ is the least upper bound of
${(x_i, r_i)}_{i \in I}$.

The inequality $d (x, y) \geq \sup_{i \in I} (d (x_i, y) - r_i)$ is
automatic: since $(x_i, r_i) \leq^{d^+} (x, 0)$, $d (x_i, x) \leq r_i$
for every $i \in I$, so
$d (x_i, y) \leq d (x_i, x) + d (x, y) \leq d (x, y) + r_i$, which
entails that $d (x_i, y) - r_i \leq d (x, y)$ for every $i \in I$.  If
the inequality were strict, there would be an $r \in \Rp$ such that
$d (x, y) > r \geq \sup_{i \in I} (d (x_i, y) - r_i)$.  We use
Proposition~\ref{prop:d:std}~(2) and note that $(x, r)$ is the
supremum of ${(x_i, r_i+r)}_{i \in I}$.  For every $i \in I$,
$d (x_i, y) - r_i \leq r$, so $(x_i, r_i+r) \leq^{d^+} (y, 0)$, and
therefore $(y, 0)$ is an upper bound of the family
${(x_i, r_i)}_{i \in I}$.  Since $(x, r)$ is the least one,
$(x, r) \leq^{d^+} (y, 0)$, that is, $d (x, y) \leq r$.  That
contradicts $d (x, y) > r$, so the inequality is an equality, showing
that $x$ is a $d$-limit of ${(x_i)}_{i \in I, \sqsubseteq}$.  \qed

\begin{thm}
  \label{thm:d:enoughcenter}
  A standard quasi-metric space $X, d$ has enough center points if and
  only if it is algebraic.
\end{thm}
\proof Assume $X, d$ is algebraic, let $\mathcal U$ be a Scott-open
subset of $\mathbf B (X, d)$ and $(y, s)$ be a formal ball in
$\mathcal U$.  By definition, $y$ is the $d$-limit of some Cauchy net
consisting of $d$-finite points, and we have seen that we could
replace that Cauchy net with a Cauchy-weightable subnet.  By
Lemma~\ref{lemma:d:limit:sup}, we can therefore express $(y, 0)$ as
the supremum of some directed family ${(x_i, r_i)}_{i \in I}$, where
each $x_i$ is $d$-finite.  By Proposition~\ref{prop:d:std:props}~(2),
$(y, s)$ is the supremum of the directed family
${(x_i, r_i+s)}_{i \in I}$.  Therefore $(x_i, r_i+s)$ is in
$\mathcal U$ for some $i \in I$.  Using Lemma~\ref{lemma:+eps}, there
is even an $\epsilon > 0$ such that $(x_i, r_i+s+\epsilon)$ is in
$\mathcal U$.  Since $(x_i, r_i+s) \leq^{d^+} (y, s)$,
$d (x_i, y) \leq (r_i+s) - s < (r_i+s)-s+\epsilon$, so $(y, s)$ is in
$B^{d^+}_{(x_i, 0), < r_i+s+\epsilon}$.  Every formal ball $(z, t)$ in
$B^{d^+}_{(x_i, 0), < r_i+s+\epsilon}$ is such that
$d (x_i, z) < r_i + s - t +\epsilon$, and that implies
$(x_i, r_i+s+\epsilon) \leq^{d^+} (z, t)$.  Since
$(x_i, r_i+s+\epsilon)$ is in $\mathcal U$, and open subsets are
upwards-closed, $(z, t)$ is in $\mathcal U$, too.  We sum up:
$(y, s) \in B^{d^+}_{(x_i, 0), < r_i+s+\epsilon} \subseteq \mathcal
U$;
since $x_i$ is $d$-finite (equivalently, a center point, see
Lemma~\ref{lemma:d:finite}), $X, d$ has enough center points.

Conversely, assume that $X, d$ has enough center points.  Fix
$x \in X$.  Let $I$ be the family of all non-empty finite sets of open
neighborhoods of $(x, 0)$ in $\mathbf B (X, d)$, ordered by set
inclusion.  To stress it, an element $i$ of $I$ is a finite set
$\{\mathcal U_1, \cdots, \mathcal U_n\}$ where $n \geq 1$ and each
$\mathcal U_i$ is a Scott-open set of formal balls containing
$(x, 0)$.  By induction on the cardinality $n$ of
$i = \{\mathcal U_1, \cdots, \mathcal U_n\}$, we build a formal ball
$(x_i, r_i)$, where $x_i$ is a center point, in such a way that
$(x, 0) \in B^{d^+}_{(x_i, 0), <r_i} \subseteq \mathcal U_1 \cap
\cdots \cap \mathcal U_n \cap \bigcap_{j \subseteq i, j\neq\emptyset,
  j\neq i} B^{d^+}_{(x_j, 0), < r_j}$.
To do so, just observe that $X, d$ has enough center points, and that,
by induction hypothesis each $B^{d^+}_{(x_j, 0), < r_j}$ is an open
neighborhood of $(x, 0)$.

If $j \subseteq i$, then
$B^{d^+}_{(x_i, 0), <r_i} \subseteq B^{d^+}_{(x_j, 0), < r_j}$.  For
every $\epsilon > 0$ such that $\epsilon \leq r_i$,
$(x_i, r_i-\epsilon)$ is in $B^{d^+}_{(x_i, 0), <r_i}$, hence in
$B^{d^+}_{(x_j, 0), < r_j}$, so $d (x_j, x_i) + r_i - \epsilon < r_j$.
As $\epsilon$ tends to $0$, $d (x_j, x_i) \leq r_j - r_i$, and this
shows that the net ${(x_i, r_i)}_{i \in I, \subseteq}$ is a monotone
net.  By construction, any upper bound $(z, t)$ of that net is in
every open neighborhood $\mathcal U$ of $(x, 0)$, hence must be above
$(x, 0)$.  By construction again, for every $i \in I$,
$(x, 0) \in B^{d^+}_{(x_i, 0), <r_i}$, so $d (x_i, x) < r_i$, which
implies $(x_i, r_i) \leq^{d^+} (x, 0)$, therefore showing that
$(x, 0)$ is an upper bound of the net.  We have seen that any upper
bound $(z, t)$ would be above $(x, 0)$, so $(x, 0)$ is the supremum of
the directed family of formal balls $(x_i, r_i)$.  Since every $x_i$
is a center point, namely, a $d$-finite point, we conclude.  \qed


As an application, recall that a quasi-metric space $X, d$ is
Smyth-complete if and only if  $\mathbf B (X, d)$ is sober in its open
ball topology.  If that is the case, then $\mathbf B (X, d)$ is a
monotone convergence space, and in particular \emph{every} open ball
of $\mathbf B (X, d), d^+$ is Scott-open.  By definition, this implies
that \emph{every} point of $X$ is a center point.

Conversely, if every point of $X$ is a center point, then the open
ball topology on $\mathbf B (X, d), d^+$ is coarser than the Scott
topology, and by Lemma~\ref{lemma:Scott=>openball}, the two topologies
coincide.  Since $\mathbf B (X, d)$ is a C-space in its open ball
topology, it is also a C-space in its Scott topology.  But the posets
that are C-spaces in their Scott topology are exactly the continuous
posets \cite[Proposition~4]{Erne:minbases}.  If $X, d$ is also
Yoneda-complete, then $\mathbf B (X, d)$ is a continuous dcpo, hence
its Scott topology is sober.  The open ball topology on
$\mathbf B (X, d)$ is then sober as well, hence $X, d$ is
Smyth-complete by the Romaguera-Valero theorem.

We have therefore obtained a proof of the following.  The equivalence
between (1) and (2) is, modulo some details, due to Ali-Akbari,
Honarii, Pourmahdian, and Rezaii \cite{AHPR:formal:ball}.  What they
state is that $X, d$ is Smyth-complete if and only if it is
Yoneda-complete and all its points are $d$-finite.  This is
equivalent, since Yoneda-complete spaces are standard, and $d$-finite
points coincide with center points in standard spaces.
\begin{prop}
  \label{lemma:Smyth-complete}
  For a quasi-metric space $X, d$, the following are equivalent:
  \begin{enumerate}
  \item $X, d$ is Smyth-complete;
  \item $X, d$ is Yoneda-complete and all its points are center
    points;
  \item $\mathbf B (X, d)$ is a dcpo, and the open ball topology on
    $\mathbf B (X, d), d^+$ coincides with the Scott topology.
  \end{enumerate}
\end{prop}

As another application, the following generalizes the fact that every
algebraic Yoneda-complete quasi-metric space is continuous.
\begin{prop}
  \label{prop:alg=>cont}
  Every (standard) algebraic quasi-metric space $X, d$ is continuous.
  When $z$ is a center point, $(z, t) \ll (y, s)$ if and only if
  $(z, t) \prec (y, s)$, if and only if $d (z, y) < t-s$.  In general,
  $(x, r) \ll (y, s)$ if and only there is a center point $z$ and some
  $t \in \Rp$ such that $(x, r) \leq^{d^+} (z, t) \prec (y, s)$.
\end{prop}
\proof For every center point $z$, for every $t \in \Rp$,
$B^{d^+}_{(z, 0), < t} = \{(y, s) \mid d (z, y) < t-s\}$ is Scott-open
by definition, and is included in $\upc (z, t)$.  Therefore
$B^{d^+}_{(z, 0), < t}$ is included in the interior of $\upc (z, t)$.
Conversely, if $(y, s)$ is in the interior of $\upc (z, t)$, then $(y,
s+1/2^n)$ is in $\upc (z, t)$ by Lemma~\ref{lemma:+eps}, so $d (z, y)
\leq t - s - 1/2^n < t-s$.

It follows that $(z, t) \ll (y, s)$ if and only if $(y, s)$ is in the
interior of $\upc (z, t)$, if and only if
$(y, s) \in B^{d^+}_{(z, 0), < t}$, if and only if $d (z, y) < t - s$.

Fix a formal ball $(y, s)$.  Since $X, d$ is algebraic, $y$ is the
$d$-limit of some Cauchy-weightable net of $d$-finite points.
That is, $(y, 0)$ is the supremum of some directed family ${(z_i,
  t_i)}_{i \in I}$ where each $z_i$ is finite (equivalently, a center
point).  Since $X, d$ is standard, $(y, s)$ is the supremum of the
directed family ${(z_i, t_i+s)}_{i \in I}$.

The family ${(z_i, t_i+s+1/2^n)}_{i \in I, n \in \nat}$ is again
directed: given $i, j \in I$ and $m, n \in \nat$, find $k \in I$ such
that $(z_i, t_i+s), (z_j, t_j+s) \leq^{d^+} (z_k, t_k+s)$, and let
$p = \min (m, n)$, then
$(z_i, t_i+s+1/2^m), (z_j, t_j+s+1/2^n) \leq^{d^+} (z_k,
t_k+s+1/2^p)$.  The upper bounds of
${(z_i, t_i+s+1/2^n)}_{i \in I, n \in \nat}$ are exactly those of
${(z_i, t_i+s)}_{i \in I}$, so
${(z_i, t_i+s+1/2^n)}_{i \in I, n \in \nat}$ admits $(y, s)$ as upper
bound.  Moreover, since $(z_i, t_i+s) \leq^{d^+} (y, s)$,
$d (z_i, y) \leq t_i < t_i+1/2^n$, hence
$(z_i, t_i+s+1/2^n) \ll (y, s)$.  This allows us to conclude that
$\mathbf B (X, d)$ is continuous, hence that $X, d$ is continuous.

Finally, for general formal balls $(x, r)$ and $(y, s)$, we show that
$(x, r) \ll (y, s)$ if and only if there is a formal ball $(z, t)$
such that $(x, r) \leq^{d^+} (z, t)$ and $d (z, y) < t-s$.  If
$(x, r) \ll (y, s)$, since $(y, s)$ is the directed supremum of a
family of formal balls $(z, t)$ ($z$ center point) way-below $(y, s)$,
$(x, r) \leq^{d^+} (z, t) \ll (y, s)$ for some such formal ball.  The
converse direction is obvious.  \qed

\section{Continuous and Lipschitz Real-Valued Maps}
\label{sec:cont-lipsch-real}

It is time we talked about morphisms.

Given $\alpha \in \Rp$, and two quasi-metric spaces $X, d$ and
$Y, \partial$, a map $f \colon X \to Y$ is \emph{$\alpha$-Lipschitz}
if and only if $\partial (f (x), f (x')) \leq \alpha \; d (x, x')$ for
all $x, x' \in X$.  (When $\alpha=0$ and $d (x, x')=+\infty$, we take
the convention that $0.+\infty = +\infty$.)  It is \emph{Lipschitz} if
and only if it is $\alpha$-Lipschitz for some $\alpha \in \Rp$.

\begin{exa}
  \label{ex:Lip:mono}
  If $X$ and $Y$ are posets, considered as quasi-metric spaces, the
  $\alpha$-Lipschitz maps for $\alpha > 0$ are the monotonic maps.  \qed
\end{exa}

We are particularly interested in the case $\alpha=1$; $1$-Lipschitz
maps are sometimes called \emph{non-expansive}.  Quasi-metric spaces
and $1$-Lipschitz maps form a category $\QMet$.

\begin{rem}
  \label{rem:expo}
  We have repeatedly taken the example of posets as specific
  quasi-metric spaces.  One can see the category of posets and
  monotonic maps as a full subcategory of $\QMet$, consisting of those
  quasi-metric spaces whose quasi-metric takes the values $0$ and
  $+\infty$ only.  It is a remarkable fact that the posets are exactly
  the exponentiable objects in $\QMet$
  \cite[Exercise~6.6.18]{JGL-topology}.
\end{rem}

Given $f \colon X \to Y$, and $\alpha \in \Rp$, we can lift $f$ to
spaces of formal balls by defining
$\mathbf B^\alpha (f) \colon \mathbf B (X, d) \to \mathbf B
(Y, \partial)$
as: $\mathbf B^\alpha (f) (x, r) = (f (x), \alpha r)$.  (In
particular, $\mathbf B^1 (f) (x, r) = (f (x), r)$.)  It is easy to
verify that $f$ is $\alpha$-Lipschitz if and only if
$\mathbf B^\alpha (f)$ is monotonic.

All Lipschitz maps are continuous with respect to the underlying open
ball topologies, but not so for $d$-Scott topologies.  Bonsangue
\emph{et al.} defined the following notion of Yoneda-continuity
\cite{BvBR:gms}.  Given a Lipschitz map $f$ as above, $f$ is
\emph{Yoneda-continuous} if and only if $f$ maps $d$-limits of Cauchy
nets in $X$ to $\partial$-limits in $Y$.  The notion applies, in more
generality, not just to Lipschitz maps, but to all uniformly
continuous maps (op.\ cit.), but we will not consider that case here.
Once again, formal balls offer a reduction of the quasi-metric concept
to a domain-theoretic notion: when $X, d$ and $Y, \partial$ are both
Yoneda-complete, and given an $\alpha$-Lipschitz map
$f \colon X \to Y$, $f$ is Yoneda-continuous if and only if
$\mathbf B^\alpha (f)$ is Scott-continuous
\cite[Proposition~7.4.38]{JGL-topology}.  That extends to standard
quasi-metric spaces:
\begin{lem}
  \label{lemma:alphaLip:cont}
  Let $X, d$ and $Y, \partial$ be two standard quasi-metric spaces.
  For every map $f \colon X \to Y$, $f$ is $\alpha$-Lipschitz
  Yoneda-continuous if and only if $\mathbf B^\alpha (f)$ is
  Scott-continuous.
\end{lem}
\proof If $f$ is $\alpha$-Lipschitz Yoneda-continuous, then
$\mathbf B^\alpha (f)$ is monotonic.  Given a formal ball $(x, r)$
that is the supremum of a directed family ${(x_i, r_i)}_{i \in I}$,
$(x, 0)$ is the supremum of the directed family
${(x_i, r_i-r)}_{i \in I}$, where $r = \inf_{i \in I} r_i$, using the
fact that $X, d$ is standard.  Hence $x$ is the $d$-limit of the
Cauchy net ${(x_i)}_{i \in I, \sqsubseteq}$, where $i \sqsubseteq j$
if and only if $(x_i, r_i) \leq^{d^+} (x_j, r_j)$.  By assumption,
$f (x)$ is the $\partial$-limit of the Cauchy net
${(f (x_i))}_{i \in I, \sqsubseteq}$.  Since $Y, \partial$ is
standard, we can use Lemma~\ref{lemma:d:limit:sup} and conclude that
$(f (x), 0)$ is the supremum of the directed family
${(f (x_i), \alpha r_i- \alpha r)}_{i \in I}$.  Using standardness
again, $(f (x), \alpha r)$ is the supremum of
${(f (x_i), \alpha r_i)}_{i \in I}$.  Therefore $\mathbf B^\alpha (f)$
is Scott-continuous.

Conversely, if $\mathbf B^\alpha (f)$ is Scott-continuous, then it is
monotonic, so $f$ is $\alpha$-Lipschitz.  For every Cauchy net
${(x_i)}_{i \in I, \sqsubseteq}$ in $X$ with a $d$-limit $x$, extract
a Cauchy-weightable subnet \cite[Lemma~7.2.8]{JGL-topology}.  We shall
show that the image of that subnet by $f$ has a $\partial$-limit,
hence the whole Cauchy net ${(f (x_i))}_{i \in I, \sqsubseteq}$ will
have the same $\partial$-limit \cite[Lemma~7.4.6]{JGL-topology}.  Note
also that the subnet we have taken still has $x$ as its $d$-limit
\cite[Exercise~7.4.7]{JGL-topology}.  (We have already used these
tricks in Section~\ref{sec:algebr-yoneda}.)  Hence, without loss of
generality, assume that ${(x_i, r_i)}_{i \in I, \sqsubseteq}$ if a
Cauchy-weighted net.  Since $x$ is a $d$-limit of our original net,
$(x, 0) = \sup_{i \in I} (x_i, r_i)$.  By Scott-continuity, $(f (x),
0) = \sup_{i \in I} (f (x_i), \alpha r_i)$.  Since $Y, \partial$ is standard,
we can use Lemma~\ref{lemma:d:limit:sup} and conclude that $f$ is the
$\partial$-limit of ${(f (x_i))}_{i \in I, \sqsubseteq}$, showing that
$f$ is Yoneda-continuous.  \qed

In turn, the above result can be simplified as follows when $Y =
\creal$.
\begin{lem}
  \label{lemma:alphaLip:cont:real}
  Let $X, d$ be a standard quasi-metric space.  A map
  $f \colon X \to \creal$ is $\alpha$-Lipschitz Yoneda-continuous from
  $X, d$ to $\creal, \dreal$ if and only if the map
  $f' \colon \mathbf B (X, d) \to \real \cup \{+\infty\}$, defined by
  $f' (x, r) = f (x) - \alpha r$, is Scott-continuous.
\end{lem}
\proof If $f$ is $\alpha$-Lipschitz Yoneda-continuous, then
$\mathbf B^\alpha (f)$ is Scott-continuous, and using the isomorphism
of Example~\ref{ex:Yoneda-complete:R}, the map
$(x, r) \mapsto (f (x) - \alpha r, -\alpha r) = (f' (x, r), -\alpha
r)$
is Scott-continuous from $\mathbf B (X, d)$ to a Scott-closed subset
of $(\real \cup \{+\infty\}) \times (-\infty, 0]$.  Taking first
components, $f'$ is Scott-continuous.

Conversely, if $f'$ is Scott-continuous, then the map
$(x, r) \mapsto (f' (x, r), -\alpha r)$ is Scott-continuous, too.  To
show this, we only need to show that the map
$(x, r) \mapsto -\alpha r$ is Scott-continuous, and that is a
consequence of Proposition~\ref{prop:d:std:props}~(3).  Using the same
isomorphism as above, this implies that $\mathbf B^\alpha (f)$ is
Scott-continuous.  \qed

\begin{fact}
  \label{fact:YQMet}
  The Yoneda-complete quasi-metric spaces, together with the
  $1$-Lipschitz Yoneda-continuous maps, form a category $\WC\QMet$.
\end{fact}

\begin{exa}
  \label{ex:Yc:Scottcont}
  If $X$ and $Y$ are posets, considered as quasi-metric spaces, the
  $\alpha$-Lipschitz Yoneda-continuous maps from $X$ to $Y$, for
  $\alpha > 0$, are the Scott-continuous maps.  \qed
\end{exa}

We shall require the following easy facts.
\begin{prop}
  \label{prop:alphaLip:props}
  Let $X, d$ be a standard quasi-metric space,
  $\alpha, \beta \in \Rp$, and $f$, $g$ be maps from $X, d$ to
  $\creal, \dreal$.
  \begin{enumerate}
  \item If $f$ is $\beta$-Lipschitz Yoneda-continuous, then $\alpha f$
    is $\alpha\beta$-Lipschitz Yoneda-continuous;
  \item If $f$ is $\alpha$-Lipschitz Yoneda-continuous and $g$ is
    $\beta$-Lipschitz Yoneda-continuous then $f+g$ is
    $(\alpha+\beta)$-Lipschitz Yoneda-continuous;
  \item If $f$, $g$ are $\alpha$-Lipschitz Yoneda-continuous, then so
    are $\min (f, g)$ and $\max (f, g)$;
  \item If ${(f_i)}_{i \in I}$ is any family of $\alpha$-Lipschitz
    Yoneda-continuous maps, then the pointwise supremum
    $\sup_{i \in I} f_i$ is also $\alpha$-Lipschitz Yoneda-continuous.
  \item If $\alpha \leq \beta$ and $f$ is $\alpha$-Lipschitz
    Yoneda-continuous then $f$ is $\beta$-Lipschitz Yoneda-continuous.
  \end{enumerate}
\end{prop}
\proof We use Lemma~\ref{lemma:alphaLip:cont:real}.  (We let the
reader check, as an exercise, that (1) and (5) hold even for
non-standard quasi-metric spaces.  But one cannot use
Lemma~\ref{lemma:alphaLip:cont:real}, then.)

(1) If $(x, r) \mapsto f (x) - \beta r$ is Scott-continuous, then $(x,
r) \mapsto \alpha f (x) - \alpha\beta r$ is, too, since multiplication
by $\alpha \in \Rp$ is Scott-continuous.  Note that it does not matter
that we multiply $\alpha$ by a positive or a negative number, but the
fact that $\alpha$ is non-negative does matter.

(2) If $(x, r) \mapsto f (x) - \alpha r$ and $(x, r) \mapsto g (x) -
\beta r$ are Scott-continuous, so is their sum $(x, r) \mapsto (f+g)
(x) - (\alpha+\beta) r$, because addition is Scott-continuous.

(3) Similarly, since $\min$ and $\max$ are Scott-continuous.

(4) The map $(x, r) \mapsto (\sup_{i \in I} f_i (x)) - \alpha r =
\sup_{i \in I} (f_i (x) - \alpha r)$ is Scott-continuous, because any
supremum of Scott-continuous maps with values in $\real \cup
\{+\infty\}$ is Scott-continuous: this is a standard exercise, and
reduces to showing that suprema commute.

(5) If $(x, r) \mapsto f (x) - \alpha r$ is Scott-continuous, then so
is $(x, r) \mapsto f (x) - \beta r$, since it arises as the sum of the
former plus the map $(x, r) \mapsto - (\beta-\alpha) r$, which is
Scott-continuous by Proposition~\ref{prop:d:std:props}~(3).  \qed

In a metric space $X, d$, there is a notion of distance to a closed
set $C$: $d (x, C) = \inf \{d (x, y) \mid y \in C\}$, and $d (x, C)=0$
if and only if $x \in C$.  One can also define the \emph{thinning} of
an open subset $U$ by $r \geq 0$ as the set of points $x$ whose
distance to the complement of $U$ is strictly larger than $r$.

We generalize the notion to all standard quasi-metric spaces as
follows.
\begin{defi}
  \label{defn:hatU}
  Given any subset $A$ of a quasi-metric space $X, d$, the largest
  open subset $V$ of $\mathbf B (X, d)$ such that $V \cap X \subseteq
  A$ is written $\widehat A$.
\end{defi}
$\widehat A$ is obtained as the union of all the open subsets $V$ of
$\mathbf B (X, d)$ such that $V \cap X \subseteq A$.  By the
definition of the $d$-Scott topology, $\widehat A \cap X$ is the
interior of $A$ in $X$.  In particular, for every $d$-Scott open
subset $U$ of $X$, $\widehat U \cap X = U$.

By Proposition~\ref{prop:d:std:props}~(4), for every $r \in \Rp$, the
map $\_ + r$ is Scott-continuous, so $(\_ + r)^{-1} (\widehat U)$ is
also open.  Hence the following defines a $d$-Scott open subset of $X$.
\begin{defi}[Thinning]
  \label{defn:thinning}
  In a standard quasi-metric space $X, d$, the \emph{thinning}
  $U^{-r}$ of the $d$-Scott open subset $U$ by $r \in \Rp$ is
  $(\_ + r)^{-1} (\widehat U) \cap X = \{x \in X \mid (x, r) \in
  \widehat U\}$.
\end{defi}



This allows us to (re)define the distance $d (x, \overline U)$ of $x$
to the complement $\overline U$ of an open subset $U$.
\begin{defi}
  \label{defn:dC}
  In a standard quasi-metric space $X, d$, define $d (x, \overline U)$
  for $x \in X$ and $U$ open in $X$ as $\sup \{r \in \Rp \mid x \in
    U^{-r}\} = \sup \{r \in \Rp \mid (x, r) \in \widehat U\}$.
\end{defi}
This has the expected properties:
\begin{lem}
  \label{lemma:dC:cont}
  In a standard quasi-metric space $X, d$, and for all $x, y \in X$
  and every $d$-Scott open subset $U$ of $X$, the following hold:
  \begin{enumerate}
  \item $d (x, \overline U)=0$ if and only if $x \not\in U$;
  \item $d (x, \overline U) \leq d (x, y) + d (y, \overline U)$;
  \item the map $d (\_, \overline U)$ is $1$-Lipschitz
    Yoneda-continuous from $X, d$ to $\creal, \dreal$.
  \end{enumerate}
\end{lem}
\proof (1) If $x \in U$, then $(x, 0) \in \widehat U$, so $(x, 1/2^n)$
is in $\widehat U$ for some $n \in \nat$, by appealing to
Lemma~\ref{lemma:+eps}.  It follows that
$d (x, \overline U) \geq 1/2^n > 0$.  Conversely, if
$d (x, \overline U) > 0$, then $(x, r) \in \widehat U$ for some
$r > 0$.  Since $(x, r) \leq^{d^+} (x, 0)$, $(x, 0)$ is also in
$\widehat U$, so $x \in U$.


(2) Assume that $d (x, \overline U) > d (x, y) + d (y, \overline U)$.
Then there is an $r \in \Rp$ such that $(x, r) \in \widehat U$ and
$r > d (x, y) + d (y, \overline U)$.  In particular, $d (x, y) < r$,
and therefore $(y, r - d (x,y))$ is a well-defined formal ball.
Moreover, $(x, r) \leq^{d^+} (y, r - d (x, y))$, by definition.  Since
$\widehat U$ is upwards-closed, $(y, r - d (x, y))$ is in
$\widehat U$, so $d (y, \overline U) \geq r - d (x, y)$:
contradiction.

(3) By (2), $f = d (\_, \overline U)$ is $1$-Lipschitz.  Hence
$f' \colon (x, r) \mapsto d (x, \overline U) - r$ is monotonic.  Let
us show that $f'$ is Scott-continuous.  We shall conclude by using
Lemma~\ref{lemma:alphaLip:cont:real}.  Let ${(x_i, r_i)}_{i \in I}$ be
a directed family of formal balls with a supremum $(x, r)$.  We must
show that
$d (x, \overline U) - r = \sup_{i \in I} d (x_i, \overline U) - r_i$.
By monotonicity, the left-hand side is larger than or equal to the
right-hand side.  Let us suppose, for the sake of contradiction, that
it is strictly larger: for some $a \in \Rp$,
$d (x, \overline U) > a > \sup_{i \in I} d (x_i, \overline U) + r -
r_i$.
By definition of $d (x, \overline U)$, there is a radius $r' \in \Rp$
such that $(x, r') \in \widehat U$, and $r' > a$.  By
Proposition~\ref{prop:d:std:props}~(2), $(x, r')$ is the supremum of
the directed family ${(x_i, r_i+r'-r)}_{i \in I}$.  Since $\widehat U$
is Scott-open, there is an $i \in I$ such that $(x_i, r_i+r'-r)$ is in
$\widehat U$.
In particular, $d (x_i, \overline U) \geq r_i + r' - r$, hence
$a > \sup_{i \in I} d (x_i, \overline U) + r - r_i \geq r'$.  This is
impossible since $r' > a$.
  \qed

The following compares $d (x, \overline U)$ with the more familiar
formula $\inf_{y \in \overline U} d (x, y)$.  We write $\dc$ for
downward closure with respect to $\leq^d$.  For a finite set $E$, $\dc
E$ is closed, hence its complement is always open.
\begin{prop}
  \label{prop:d:Hoare}
  Let $X, d$ be a standard quasi-metric space, $x \in X$, and $U$ be a
  $d$-Scott open subset of $X$.  Then
  $d (x, \overline U) \leq \inf_{y \in \overline U} d (x, y)$, with
  equality if $\overline U = \dc E$ for some finite set $E$, or if $x$
  is a center point.
\end{prop}
In particular, $d (x, \overline U) =\inf_{y \in \overline U} d (x,
y)$ in all metric spaces, and in all Smyth-complete quasi-metric
spaces, since all points are center points in those situations.

\proof For every $r \in \Rp$, for every $y \in X$, if
$d (x, y) \leq r$ then $(x, r) \leq^{d^+} (y, 0)$.  Hence if
$d (x, y) \leq r$ and $(x, r) \in \widehat U$ then
$(y, 0) \in \widehat U$, hence $y \in U$.  By contraposition, if
$(x, r) \in \widehat U$ and $y \in \overline U$, then $d (x, y) > r$.
Taking infima over $y$ and suprema over $r$, we obtain
$d (x, \overline U) \leq \inf_{y \in \overline U} d (x, y)$.

If $\overline U = \dc E$ for some finite set
$E = \{y_1, y_2, \cdots, y_n\}$, then the downward closure of $E$ in
$\mathbf B (X, d)$, which we shall write as $\dc_{\mathbf B} E$ to
avoid any confusion, is the closure of $E$ in $\mathbf B (X, d)$, and
its intersection with $X$ is $\dc E$.  It follows that $\widehat U$ is
the complement of $\dc_{\mathbf B} E$.  For every
$r < \min_{1\leq i\leq n} d (x, y_i)$, for every $i$, $1\leq i\leq n$,
$(x, r)$ is not below $(y_i, 0)$, since that would imply
$d (x, y_i) \leq r$.  This means that $(x, r)$ is not in
$\dc_{\mathbf B} E$, hence is in its complement, $\widehat U$.  By
definition, it follows that $d (x, \overline U) \geq r$.  As $r$ is
arbitrary, $d (x, \overline U) \geq \min_{1\leq i\leq n} d (x, y_i)
\geq \inf_{y \in \overline U} d (x, y)$.
%

For the second part, we no longer assume $\overline U = \dc E$, but we
assume that $x$ is a center point.  We know that
$d (x, \overline U) \leq \inf_{y \in \overline U} d (x, y)$, and we
assume that the inequality is strict: there are two real numbers
$s, t \in \Rp$ such that
$d (x, \overline U) < s < t \leq \inf_{y \in \overline U} d (x, y)$.
The rightmost inequality states that every $y \in \overline U$ is such
that $d (x, y) \geq t$, hence, by contraposition, that
$B^d_{x, < t} \subseteq U$.  Since $x$ is a center point,
$B^d_{(x, s), < t-s}$ is open, by Lemma~\ref{lemma:d:finite}~(2).  The
intersection of $B^d_{(x, s), < t-s}$ with $X$ is $B^d_{x, < t}$,
which is included in $U$, so
$B^d_{(x, s), < t-s} \subseteq \widehat U$.  In particular,
$(x, s) \in \widehat U$, so $d (x, \overline U) \geq s$,
contradiction.  \qed

\begin{figure}
  \centering
  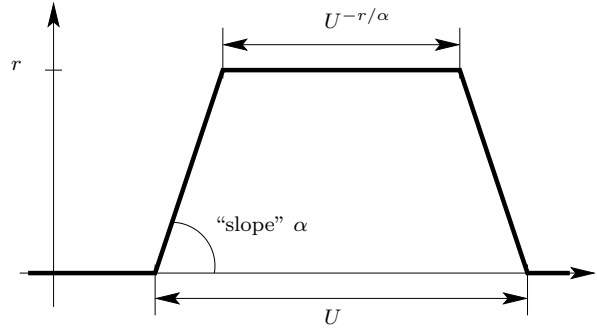
  \caption{An $\alpha$-Lipschitz Yoneda-continuous map approximating
    $r\chi_U$}
  \label{fig:rdalpha}
\end{figure}
Write $\chi_U$ for the characteristic function of the open subset $U$.
We compare functions, and take suprema of functions, pointwise.
The map $\min (r, \alpha d (\_, \allowbreak \overline U))$ studied
below is probably best understood through a picture: see Figure~\ref{fig:rdalpha}.
\begin{prop}
  \label{prop:chiU:Lip}
  Let $X, d$ be a standard quasi-metric space.  For all
  $\alpha, r \in \Rp$,
  $\min (r, \alpha d (\_, \allowbreak \overline U))$ is
  an $\alpha$-Lipschitz Yoneda-continuous map from $X, d$ to
  $\creal, \dreal$, and is less than or equal to $r \chi_U$.  Moreover
  the family ${(\min (r, \alpha d (\_, \overline U)))}_{\alpha > 0}$
  is a chain, and its supremum is $r \chi_U$.
\end{prop}
\proof The function $\min (r, \alpha d (\_, \overline U))$ is
$\alpha$-Lipschitz Yoneda-continuous by Lemma~\ref{lemma:dC:cont}~(3)
and Proposition~\ref{prop:alphaLip:props}.


For every $x \in U$,
$\min (r, \alpha d (x, \overline U)) \leq r = r \chi_U (x)$.  For
every $x \in \overline U$, we use Lemma~\ref{lemma:dC:cont}~(1) to
conclude that
$\min (r, \alpha d (x, \overline U)) = \min (r, 0) = 0 \leq \chi_U
(x)$.

If $\alpha \leq \alpha'$, then clearly $\min (r, \alpha d (x,
\overline U)) \leq \min (r, \alpha' d (x, \overline U))$, so the
family is a chain.

To show the final claim, take any $x \in U$.  By
Lemma~\ref{lemma:dC:cont}~(1), $d (x, \overline U)$ is non-zero, so
$\alpha d (x, \overline U) \geq r$ for $\alpha$ large enough.  Then
$\min d (r, \alpha d (x, \overline U)) = r = r \chi_U (x)$.  \qed

\begin{prop}
  \label{prop:chiU:largestLip}
  Let $X, d$ be a standard quasi-metric space.  For all
  $\alpha, r \in \Rp$,
  $\min (r, \alpha d (\_, \allowbreak \overline U))$ is the largest
  $\alpha$-Lipschitz Yoneda-continuous map from $X, d$ to
  $\creal, \dreal$ that is less than or equal to $r \chi_U$.
\end{prop}
\proof The claim is clear if $r=0$, so let us assume $r > 0$.  Let $f$
be $\alpha$-Lipschitz Yoneda-continuous from $X, d$ to
$\creal, \dreal$, and assume $f \leq r \chi_U$.  Recall that the map
$f' \colon \mathbf B (X, d) \to \real \cup \{+\infty\}$ defined by
$f' (x, s) = f (x) - \alpha s$ is Scott-continuous
(Lemma~\ref{lemma:alphaLip:cont:real}).

Consider the open subset $V = {f'}^{-1} (0, +\infty]$.  For every open
ball $(x, 0)$ with radius $0$ in $V$, $f' (x, 0) = f (x) > 0$.  Since
$f \leq r \chi_U$, $x$ is then in $U$.  This shows that
$V \cap X \subseteq U$, hence $V \subseteq \widehat U$.  Said in
another way, for every open ball $(x, s)$ such that
$f (x) - \alpha s > 0$, $(x, s)$ is in $\widehat U$.  Therefore, for
every $x \in X$, every $s \in \Rp$ such that $f (x) > \alpha s$ is
less than or equal to $d (x, \overline U)$.  Taking suprema over $s$,
$f (x) \leq \alpha d (x, \overline U)$.  Hence
$f \leq \alpha d (\_, \overline U)$, and we conclude since
$f \leq r \chi_U \leq r$.  \qed

A map $f \colon X \to \creal$ that is continuous when $\creal$ is
equipped with its Scott topology is classically known as a \emph{lower
  semicontinuous} function from $X$ to $\creal$.

We finally obtain the following result, which shows that we can
approximate any $\creal$-valued lower semicontinuous map, as closely
as we wish, by $\alpha$-Lipschitz Yoneda-continuous maps, as $\alpha$
tends to $+\infty$.
\begin{defi}
  \label{defn:falpha}
  For every lower semicontinuous map $f$ from a standard quasi-metric
  space $X, d$ to $\creal$, for every $\alpha \in \Rp$, let
  $f^{(\alpha)}$ be the largest $\alpha$-Lipschitz Yoneda-continuous
  map from $X, d$ to $\creal, \dreal$ below $f$.
\end{defi}
This exists, as the pointwise supremum of all $\alpha$-Lipschitz
Yoneda-continuous maps below $f$
(Proposition~\ref{prop:alphaLip:props}).
Proposition~\ref{prop:chiU:largestLip} can be recast as follows.
\begin{fact}
  \label{fact:chiU:largestLip}
  For every $r \in \Rp$, for every $d$-Scott open subset of $X$, $(r
  \chi_U)^{(\alpha)} = \min (r, \alpha d (\_, \allowbreak \overline U))$.
\end{fact}

\begin{thm}
  \label{thm:alphaLip:approx}
  Let $X, d$ be a standard quasi-metric space.  For every lower
  semicontinuous map $f \colon X \to \creal$, the family
  ${(f^{(\alpha)})}_{\alpha \in \Rp}$ is a chain, and
  $\sup_{\alpha \in \Rp} f^{(\alpha)} = f$.
%
\end{thm}
\proof The family is non-empty, since for example the constant $0$ map
is in it.  If $\alpha \leq \beta$, then $f^{(\alpha)}$ is
$\beta$-Lipschitz Yoneda-continuous by
Proposition~\ref{prop:alphaLip:props}~(5), and since $f^{(\beta)}$ is
largest, $f^{(\alpha)} \leq f^{(\beta)}$.  Hence the family is a
chain, and is in particular directed.

Clearly, $\sup_{\alpha \in \Rp} f^{(\alpha)} \leq f$.  If the
inequality were strict, there would be a point $x \in X$ and two real
numbers $r, s\in \Rp$ such that, for every $\alpha \in \Rp$,
$f^{(\alpha)} (x) \leq r < s < f (x)$.  Let $U$ be the open set
$f^{-1} (s, +\infty]$.  Then $s \chi_U \leq f$, so
$(s \chi_U)^{(\alpha)} \leq f^{(\alpha)}$ for every $\alpha \in \Rp$.
Using Fact~\ref{fact:chiU:largestLip} and
Proposition~\ref{prop:chiU:Lip},
$\sup_\alpha (s\chi_U)^{(\alpha)} = s \chi_U$, so
$s \chi_U (x) \leq \sup_\alpha f^{(\alpha)} (x) \leq r$.  This is
impossible, since $x \in U$.  \qed

Our intended application of this result is the following.  Given a
topological space $X$, let $\mathcal L X$ be the set of lower
semicontinuous maps from $X$ to $\creal$.  When $X, d$ is a standard
quasi-metric space, let $\mathcal L_1 X$ be the set of 1-Lipschitz
Yoneda-continuous maps from $X, d$ to $\creal, \dreal$.  A
\emph{prevision} on a topological space $X$ is a Scott-continuous map
$F \colon \mathcal L X \to \creal$ such that
$F (\alpha h) = \alpha F (h)$ for all $\alpha \in \Rp$,
$h \in \mathcal L X$.  Various refinements of the notion yield
semantic models for mixed probabilistic and non-deterministic choice
(see \cite{Gou-csl07}).

Define the following variant of the Hutchinson-Kantorovitch metric,
itself inspired from \cite{Gou-fossacs08b}.  The only difference is
that $h$ is not restricted to be $1$-Lipschitz, but $1$-Lipschitz and
Yoneda-continuous:
\[
d_{\mathrm H} (F, F') = \sup \{\dreal (F (h), F' (h)) \mid h \in
\mathcal L_1 X\}.
\]
A complete study of that quasi-metric is out of scope of this paper,
but showing that it is a quasi-metric at all requires
Theorem~\ref{thm:alphaLip:approx}.  It satisfies the triangular
inequality since $\dreal$ does, and the challenge is to show that
$d_{\mathrm H} (F, F') = d_{\mathrm H} (F', F) = 0$ if and only if
$F=F'$.  We show the more general claim that
$d_{\mathrm H} (F, F') = 0$ if and only if $F \leq F'$, i.e., if and
only if $F (h) \leq F' (h)$ for every $h \in \mathcal L X$.  The if
direction is obvious, while in the only if direction, we have
$\dreal (F (1/\alpha\; h^{(\alpha)}), F' (1/\alpha\; h^{(\alpha)})) = 0$
hence $F (h^{(\alpha)}) \leq F' (h^{(\alpha)})$ for every
$\alpha > 0$; using Theorem~\ref{thm:alphaLip:approx} and the
Scott-continuity of $F$ and $F'$, $F (h) \leq F' (h)$.

\section{Continuous and Algebraic Quasi-Metric Spaces}
\label{sec:cont-algebr-quasi}

Let us return for a moment to continuous and algebraic quasi-metric
spaces.  It is well-known that the continuous dcpos are exactly the
retracts of algebraic dcpos.  We shall prove a similar result for
Yoneda-complete quasi-metric spaces.

Every quasi-metric space $X, d$ has a so-called
\emph{Yoneda-completion} $\mathbf Y (X, d)$ \cite{BvBR:gms}.  It is
built as a certain subspace of the space of all $1$-Lipschitz maps
from $X, d$ to $\creal, \dreal$, and comes with an isometric embedding
$\eta^{\mathbf Y}_X \colon X \to \mathbf Y (X, d)$, which maps $x$ to
$d (\_, x)$.  That allows us to see $X, d$ as a sub-quasi-metric space
of $\mathbf Y (X, d)$.  $\mathbf Y (X, d)$ is algebraic, with $X$ as
set of $d$-finite elements.  It has the following universal property:
for every Yoneda-complete space $Y, \partial$, every
$\alpha$-Lipschitz map $f \colon X, d \to Y, \partial$ has a unique
$\alpha$-Lipschitz Yoneda-continuous extension to $\mathbf Y (X, d)$.
This holds, more generally, if one replaces `$\alpha$-Lipschitz' by
`uniformly continuous' \cite{BvBR:gms}.

There is another notion of completion, the \emph{formal ball
  completion} $\FB (X, d)$ of $X, d$: that one is due to Vickers
\cite{Vickers:completion:gms:I}; we rely on the presentation of
\cite[Section~7.5]{JGL-topology}.  This rests on a familiar
domain-theoretic construction known as the \emph{rounded ideal
  completion} $\RIdl (B, \prec)$ of an abstract basis $B, \prec$.  The
elements of $\RIdl (B, \prec)$ are the \emph{rounded ideals} of $B$,
namely subsets $D$ of $B$ that are $\prec$-directed (for every finite
subset $\{x_1, \cdots, x_n\}$ of $D$, there is an $x \in D$ such that
$x_i \prec x$ for every $i$) and $\prec$-downwards-closed (if
$y \in D$ and $x \prec y$ then $x \in D$).  $\RIdl (B, \prec)$ is
always a continuous dcpo under inclusion.  $B$ embeds into $\RIdl (B)$
through the map $b \mapsto \Dc b = \{b' \in B \mid b' \prec b\}$.
Modulo that embedding, $B$ is a basis of $\RIdl (B, \prec)$, and the
way-below relation of $\RIdl (B, \prec)$ restricted to $B$ is exactly
$\prec$.  Moreover, when $B$ is a C-space, hence an abstract basis
with $b \prec b'$ if and only if $b'$ is in the interior of $\upc b$,
then $\RIdl (B, \prec)$ is the sobrification of $B$, a result due to
J. Lawson \cite{Lawson:RIdl}.

This applies notably to $\mathbf B (X, d)$, which is a C-space in its
open ball topology, hence an abstract basis with $\prec$ defined by
$(x, r) \prec (y, s)$ if and only if $d (x, y) < r-s$.  Given a
rounded ideal $D$ of $\mathbf B (X, d)$, let its \emph{aperture}
$\alpha (D)$ be $\inf \{r \mid (x, r) \in D\}$.  Then $X$ embeds into
$\RIdl (\mathbf B (X, d), \prec)$ through the map
$x \mapsto \Dc (x, 0)$.  Clearly, $\Dc (x, 0)$ has aperture $0$.  We
now define $\FB (X, d)$ as the subset of
$\RIdl (\mathbf B (X, d), \prec)$ consisting of the rounded ideals
with aperture zero.  This comes with a quasi-metric $d^+_\Hoare$,
defined by
$d^+_\Hoare (D, D') = \sup_{b \in D} \inf_{b' \in D'} d^+ (b, b')$,
which makes the embedding $x \mapsto \Dc (x, 0)$ an isometric
embedding---namely,
$d^+_\Hoare (\Dc (x, 0), \Dc (x', 0)) = d (x, x')$.

The key property of that construction is that
$\mathbf B (\FB (X, d), d^+_\Hoare)$ is isomorphic, as a poset, to the
rounded ideal completion
$\RIdl (\mathbf B (X, d), \allowbreak \prec)$, through the map
$\sigma_X \colon (D, s) \mapsto D+s = \{(x, r+s) \mid (x, r) \in D\}$
\cite[Proposition~7.5.4]{JGL-topology}.  It follows immediately that
$\mathbf B (\FB (X, d), d^+_\Hoare)$ is a continuous dcpo, so
$\FB (X, d)$ is continuous Yoneda-complete (Theorem~\ref{thm:cont}).
In fact, it is algebraic Yoneda-complete, and the elements
$\Dc (x, 0)$ are its $d$-finite elements
\cite[Exercise~7.5.11]{JGL-topology}.

$\FB (X, d)$ has the same universal property as $\mathbf Y (X, d)$: for
every Yoneda-complete space $Y, \partial$, every $\alpha$-Lipschitz
map $f \colon X, d \to Y, \partial$ has a unique $\alpha$-Lipschitz
Yoneda-continuous extension to $\FB (X, d)$.  This holds, more
generally, if one replaces `$\alpha$-Lipschitz' by `uniformly
continuous' \cite[Proposition~7.5.22, Fact~7.5.23]{JGL-topology}.

By categorical generalities, it follows that $\FB (X, d)$ and
$\mathbf Y (X, d)$ are naturally isomorphic in $\WC\QMet$
\cite[Exercise~7.5.30]{JGL-topology}.  If we restrict ourselves to the
full subcategories of $\QMet$ and $\WC\QMet$ consisting of metric
spaces, the same argument shows that $\FB (X, d)$ and
$\mathbf Y (X, d)$ are naturally isomorphic to the usual Cauchy
completion of $X, d$.  This requires one to check that $d^+_\Hoare$ is
a metric whenever $d$ is, and that is indeed the case
\cite[Lemma~7.5.17]{JGL-topology}.

If instead we restrict ourselves to the full subcategory $\Ord$ of
$\QMet$ consisting of posets and monotonic maps, and to the full
subcategory $\Dcpo$ of $\WC\QMet$ consisting of dcpos and
Scott-continuous maps, yet the same argument, using the fact that the
ideal completion $\Idl (B)$ of a poset $B$ is the free dcpo over $B$,
leads to the conclusion that $\FB (X, d_\leq)$ and
$\mathbf Y (X, d_\leq)$ are naturally isomorphic to the ideal
completion of the poset $X$ \cite[Exercise~7.5.25]{JGL-topology}.

All that is known.  We would like to show that there is little more to
do to obtain a few interesting new results.

In any category $\mathbf C$, there is a notion of
($\mathbf C$-)\emph{retraction} of an object $Y$ onto an object $X$: a
pair of morphisms $\mathfrak r \colon Y \to X$ and
$\mathfrak s \colon X \to Y$ such that
$\mathfrak r \circ \mathfrak s = \mathrm{id}_X$.  $X$ is a
\emph{retract} of $Y$, $\mathfrak s$ is the \emph{section} map, and
$\mathfrak r$ itself is sometimes called a \emph{retraction} map.
Colloquially, we shall call a $\Dcpo$-retract a Scott-continuous
retract, and a $\WC\QMet$-retract a $1$-Lipschitz Yoneda-continuous
retract.

Recall the construction
$\mathbf B^1 (f) \colon (x, r) \in \mathbf B (X, d) \mapsto (f (x), r)
\in \mathbf B (Y, \partial)$,
for each $1$-Lipschitz map $f \colon X, d \to Y, \partial$.  This
defines a functor $\mathbf B^1$ from $\QMet$ to $\Ord$.  Recall that,
for a $1$-Lipschitz map $f$ between Yoneda-complete spaces, $f$ is
Yoneda-continuous if and only if $\mathbf B^1 (f)$ is Scott-continuous
\cite[Proposition~7.4.38]{JGL-topology}.  It follows that
$\mathbf B^1$ also defines a functor from $\WC\QMet$ to $\Dcpo$.

In particular, given a retraction
$\xymatrix{Y, \partial \ar@<1ex>[r]^{\mathfrak r} & X, d
  \ar@<1ex>[l]^{\mathfrak s}}$
in $\WC\QMet$, we obtain a retraction
$\xymatrix{\mathbf B (Y, \partial) \ar@<1ex>[r]^{\mathbf B^1
    (\mathfrak r)} & \mathbf B (X, d) \ar@<1ex>[l]^{\mathbf B^1
    (\mathfrak s)}}$
in $\Dcpo$.  The Scott-continuous retracts of continuous dcpos are
continuous dcpos, hence the following is obvious, in the light of
Theorem~\ref{thm:cont}.
\begin{prop}
  \label{prop:proj:cont}
  Any $1$-Lipschitz Yoneda-continuous retract (i.e., any retract in
  $\WC\QMet$) of a continuous Yoneda-complete quasi-metric space is
  itself continuous Yoneda-complete.
\end{prop}

Proposition~\ref{prop:proj:cont} is due to P. Waszkiewicz, in the more
general setting of domains over a Girard quantale
\cite[Theorem~3.3]{Waskiewicz:domains:girard:quantales}.

Proposition~\ref{prop:proj:cont}, together with
Proposition~\ref{prop:alg=>cont}, implies that every $1$-Lipschitz
Yoneda-continuous retract of a (standard) algebraic quasi-metric space
is continuous.

We now proceed to show that every continuous Yoneda-complete space
$X, d$ is a retract of an algebraic Yoneda-complete space, and that
will be $\FB (X, d), d^+_\Hoare$.  That seems to be new.
\begin{lem}
  \label{lemma:ll:prec}
  Let $X, d$ be a quasi-metric space.  If $(x, r) \ll (y, s)$ then
  $(x, r) \prec (y, s)$.
\end{lem}
\proof By Lemma~\ref{lemma:+eps}, $(x, r) \leq^{d^+} (y, s+1/2^n)$ for
some $n \in \nat$.  Hence $d (x, y) \leq r - s - 1/2^n < r-s$.  \qed

\begin{lem}
  \label{lemma:dd}
  Let $X, d$ be a quasi-metric space with a continuous poset of formal
  balls.  For every formal ball $(x, r)$, $\ddarrow (x, r)$ is an
  element of $\RIdl (\mathbf B (X, d), \allowbreak \prec)$.  If $r=0$,
  and $X, d$ is standard, then it is an element of $\FB (X, d)$.
\end{lem}
The notation $\ddarrow (x, r)$ stands for $\{(y, s) \in \mathbf B (X,
d) \mid (y, s) \ll (x, r)\}$.  This should not be confused with $\Dc
(x, r) = \{(y, s) \in \mathbf B (X, d) \mid (y, s) \prec (x, r)\}$.

\proof Since $\mathbf B (X, d)$ is a continuous poset,
$\ddarrow (x, r)$ is $\ll$-directed, hence $\prec$-directed by
Lemma~\ref{lemma:ll:prec}.  For every $(y, s) \in \ddarrow (x, r)$,
for every $(z, t) \prec (y, s)$, we have
$(z, t) \leq^{d^+} (y, s) \ll (x, r)$, so $(z, t) \in \ddarrow (x, r)$.
Therefore $\ddarrow (x, r)$ is $\prec$-downwards-closed, hence in
$\RIdl (\mathbf B (X, d), \allowbreak \prec)$.

If $r=0$ and $X, d$ is standard, since $(x, 0)$ is the supremum of
$\ddarrow (x, 0)$, the infimum of the radii of formal balls in
$\ddarrow (x, 0)$ is equal to $0$, by
Proposition~\ref{prop:d:std:props}~(1).  Hence the aperture of
$\ddarrow (x, 0)$ is $0$, whence $\ddarrow (x, 0) \in \FB (X, d)$.
\qed

\begin{lem}
  \label{lemma:dd+r}
  Let $X, d$ be a quasi-metric space, and assume that the way-below
  relation $\ll$ on $\mathbf B (X, d)$ is standard.  For every formal
  ball $(x, r)$, $\ddarrow (x, r) = \ddarrow (x, 0) + r$.
\end{lem}
\proof If $(y, s) \ll (x, r)$, then in particular
$(y, s) \leq^{d^+} (x, r)$, hence $s \geq r$.  Using
Definition~\ref{defn:std} (or Lemma~\ref{lemma:ll:std:half}) with
$a=r$, we obtain $(y, s-r) \ll (x, 0)$, and that exhibits $(y, s)$ as
an element of $\ddarrow (x, 0) + r$.  Conversely, any element
$(y, s+r)$ of $\ddarrow (x, 0) + r$, that is with $(y, s) \ll (x, 0)$,
satisfies $(y, s+r) \ll (x, r)$ since $\ll$ is standard.  \qed

Instead of embedding $X$ into $\FB (X, d)$ through $\eta_{\FB} (x) =
\Dc (x, 0)$, we consider $\eta'_{\FB} (x) = \ddarrow (x, 0)$.
\begin{lem}
  \label{lemma:eta'}
  Let $X, d$ be a continuous
  quasi-metric space.  The map
  $\eta'_{\FB} \colon X, d \to \FB (X, d), d^+_\Hoare$ is
  $1$-Lipschitz Yoneda-continuous, and $\mathbf B^1 (\eta'_{\FB})$
  maps $(x, r)$ to $\sigma_X^{-1} (\ddarrow (x, r))$.
\end{lem}
\proof Lemma~\ref{lemma:dd} enables us to claim that $\eta'_{\FB}$
takes its values in $\FB (X, d)$.  Now consider the map
$f \colon \mathbf B (X, d) \to \RIdl (\mathbf B (X, d), \allowbreak
\prec)$
defined by $f ((x, r)) = \ddarrow (x, r)$.  Since $\mathbf B (X, d)$
is a continuous poset, this is a Scott-continuous map.  Using
Lemma~\ref{lemma:dd+r}, we obtain that
$f ((x, r)) = \ddarrow (x, 0) + r = \sigma_X (\ddarrow (x, 0), r) =
\sigma_X (\mathbf B^1 (\eta'_{\FB}) (x, r))$.
Since $\sigma_X$ is an isomorphism, this suffices to show that
$\mathbf B^1 (\eta'_{\FB})$ is Scott-continuous, and we have seen that
this is equivalent to the fact that $\eta'_{\FB}$ is $1$-Lipschitz
Yoneda-continuous, in the case of standard quasi-metric spaces.
Recall from Definition~\ref{defn:cont} that all continuous
quasi-metric spaces are standard.  \qed

In the converse direction, given a Yoneda-complete quasi-metric space
$X, d$, let $d\text{-lim}$ map every
$D = {(x_i, r_i)}_{i \in I, \sqsubseteq} \in \FB (X, d)$ to the
$d$-limit of the Cauchy-weightable net
${(x_i)}_{i \in I, \sqsubseteq}$.  Let also $\sup$ denote the supremum
map from $\RIdl (\mathbf B (X, d), \prec)$ to $\mathbf B (X, d)$.
\begin{lem}
  \label{lem:dlim}
  Let $X, d$ be a continuous Yoneda-complete quasi-metric space.  The
  map $d\text{-lim} \colon \FB (X, d), d^+_\Hoare \to X, d$ is
  $1$-Lipschitz Yoneda-continuous, and
  $\mathbf B^1 (d\text{-lim}) = \sup \circ \sigma_X$.
\end{lem}
\proof Every element of $\RIdl (\mathbf B (X, d), \prec)$ can be
written uniquely as $D+r$, where $D \in \FB (X, d)$ and
$r = \alpha (D)$.  We know that $\sup (D+r) = (d\text{-lim} (D), r)$,
and this shows that
$\sup \circ \sigma_X = \mathbf B^1 (d\text{-lim})$.  The claim follows
immediately.  \qed

\begin{lem}
  \label{lem:dlim:eta'}
  Let $X, d$ be a Yoneda-complete quasi-metric space, and assume that
  the way-below relation $\ll$ on $\mathbf B (X, d)$ is standard.
  Then $d\text{-lim} \circ \eta'_{\FB}$ is the identity on $X$, and
  $\eta'_{\FB} \circ d{\text{-lim}} \leq \identity {\FB (X, d)}$.
\end{lem}
\proof We note that
$\mathbf B^1 (d\text{-lim} \circ \eta'_{\FB}) = \mathbf B^1
(d\text{-lim}) \circ \mathbf B^1 (\eta'_{\FB}) = (\sup \circ \sigma_X)
\circ (\sigma_X^{-1} \circ f)$,
where $f \colon (x, r) \mapsto \ddarrow (x, r)$.  Clearly,
$\sup \circ f = \identity {\mathbf B (X, d)}$.  Hence
$\mathbf B^1 (d\text{-lim} \circ \eta'_{\FB}) = \identity {\mathbf B
  (X, d)}$.
By applying each side of the equation to $(x, 0)$, we obtain that
$d\text{-lim} (\eta'_{\FB} (x)) = x$.

For the second claim, we consider any rounded ideal $D$ of
$\mathbf B (X, d), \prec$.  Let $x = d\text{-lim} (D)$, so that
$(x, 0) = \sup D$.  For every
$(y, s) \in \eta'_{\FB} (x) = \ddarrow (x, 0)$, use interpolation to
find $(z, t)$ such that $(y, s) \ll (z, t) \ll (x, 0) = \sup D$, so
$(y, s)$ is way-below some element $(y', s')$ of $D$.  By
Lemma~\ref{lemma:ll:prec}, $(y, s) \prec (y', s')$, so $(y, s) \in D$
since $D$ is $\prec$-downwards-closed.  Therefore
$\eta'_{\FB} (d\text{-lim} (D)) \subseteq D$.  \qed

The second part of Lemma~\ref{lem:dlim:eta'} shows that $d\text{-lim}$
and $\eta'_{\FB}$ not only define a retraction, but an
\emph{embedding-projection pair} (a concept that has meaning in any
order-enriched category; here, $\WC\QMet$).  A retract defined this
way is called a \emph{projection}.

Putting all this together, we obtain:
\begin{prop}
  \label{prop:FB:retr}
  Every continuous Yoneda-complete quasi-metric space $X, d$ is a
  projection of the algebraic Yoneda-complete quasi-metric space
  $\FB (X, d), d^+_\Hoare$ through the pair $\eta'_{\FB}$,
  $d\text{-lim}$.
\end{prop}


\begin{thm}
  \label{thm:retract:cont}
  The continuous Yoneda-complete quasi-metric spaces are exactly the
  $1$-Lipschitz Yoneda-continuous retracts (resp., projections) of
  algebraic Yoneda-complete quasi-metric spaces.
\end{thm}

\begin{rem}
  \label{rem:retr:cont}
  Theorem~\ref{thm:retract:cont} is very similar to the well-known
  result that the continuous dcpos are exactly the Scott-continuous
  retracts (resp., projections) of algebraic dcpos, and our proof is
  also very similar.  The main difference is our use of a rounded
  ideal completion instead of an ideal completion.  In fact,
  Theorem~\ref{thm:retract:cont} includes that domain-theoretic result
  as a special case.  Notably, if $X$ is a dcpo, we can consider it as
  a Yoneda-complete quasi-metric space $X, d_\leq$.  Then
  $\FB (X, d_\leq)$ is easily seen to be exactly the ideal completion
  of $X$, and ${d_\leq}^+_\Hoare$ is exactly $d_\subseteq$.
\end{rem}

\section{Quasi-Ideal Models}
\label{sec:quasi-ideal-models}

Keye Martin introduced the notion of an \emph{ideal domain}
\cite{Martin:idealmodels}, namely dcpos where each non-finite element
is maximal.  All such domains are automatically algebraic, and
first-countable.  If we agree that a \emph{model} of a space $X$ is a
dcpo in which $X$ embeds as its space of maximal elements, Martin also
showed that every space $X$ that has an $\omega$-continuous model has
an ideal model; and that the metrizable spaces that have an ideal
model are exactly the completely metrizable spaces.

By definition, a space that has a model must be $T_1$.  It is tempting
to try and generalize the notion of model to $T_0$ spaces, say as a
dcpo in which $X$ embeds as an upwards-closed subspace.
To generalize ideal models, we shall require $X$ to embed as the set
of non-finite elements in a quasi-ideal domain, defined as follows.
\begin{defi}[Quasi-Ideal Domain]
  \label{defn:qideal}
  A \emph{quasi-ideal domain} is an algebraic domain in which every
  element below a finite element is itself finite.  An
  \emph{$\omega$-quasi-ideal domain} is a quasi-ideal domain that has
  only countably many finite elements.
\end{defi}
Ideal domains are clearly quasi-ideal.  In a quasi-ideal domain, we
shall call \emph{limit elements} those points that are not finite.  A
quasi-ideal domain is organized as two non-mixing layers: a layer of
finite elements, all below a second layer of limit elements.
\begin{exa}
  \label{exa:qideal}
  For a quasi-ideal domain that is not ideal, consider $\pow (A)$
  under inclusion, for any infinite set $A$.  This is an
  $\omega$-quasi-ideal domain if and only if $A$ is countable.  \qed
\end{exa}
\begin{exa}
  \label{exa:qideal:no}
  Any quasi-ideal domain is isomorphic to the ideal completion of its
  poset of finite elements, because that is the case for all algebraic
  domains.  Conversely, given a poset $B$, its ideal completion
  $\Idl (B)$ is algebraic, but almost never quasi-ideal.  For example,
  $\Idl (\Rp)$ consists of finite elements of the form $[0, a]$, and
  limit elements $[0, b)$, $b > 0$.  The order is inclusion, and they
  are deeply interleaved.  In general, we can show that $\Idl (B)$ is
  a quasi-ideal domain if and only if, for every directed family $D$
  in $B$ with no largest element, $D$ has no upper bound in $B$.
  Indeed, if $D$ has no largest element, then $\dc D$ is a limit
  element in $\Idl (B)$, and if $D$ has an upper bound $x$, this limit
  element is below the finite element $\dc x$, showing that $\Idl (B)$
  is not quasi-ideal.  Conversely, if the directed families in $B$
  with no largest element have no upper bound, then for every ideal
  $D$ in $\Idl (B)$ that is a limit element, there is no finite
  element $\dc x$, $x \in B$, such that $D \subseteq \dc x$, so that
  $\Idl (B)$ is a quasi-ideal domain.  That observation simplifies to
  the following when $B$ is a dcpo: for a dcpo $B$, $\Idl (B)$ is a
  quasi-ideal domain if and only if $B$ has the ascending chain
  condition.  \qed
\end{exa}

An ($\omega$-)\emph{quasi-ideal model} of a topological space $X$ is
an ($\omega$-)quasi-ideal domain, seen as a topological space with the
Scott topology, whose subspace of limit elements is homeomorphic to
$X$ ---~in short, an ($\omega$-)quasi-ideal domain in which $X$ embeds
as its subspace of limit elements.


De Brecht showed \cite[Theorem~53]{deBrecht:qPolish} that the
quasi-Polish spaces are exactly the spaces that embed as the
non-finite elements of some $\omega$-algebraic (equivalently here,
$\omega$-continuous) domain.  One consequence of our results below
will be a strengthening of one direction of that theorem, namely that
all quasi-Polish spaces have an $\omega$-quasi-ideal model.  We also
believe that the proof is simpler.

For now, our goal will be slightly different: to show that every
continuous Yoneda-complete quasi-metric space has a quasi-ideal model;
but the technique will be the same.  The basic construction is
inspired by what Martin did \cite{Martin:idealmodels}.

Given a continuous Yoneda-complete quasi-metric space $X, d$,
$\mathbf B (X, d)$ has a basis.  We need slightly less than that.
\begin{defi}
  \label{defn:locbasis}
  Given a quasi-metric space $X, d$ with a continuous poset of formal
  balls, a \emph{local basis} of $\mathbf B (X, d)$ is a subset
  $\mathcal B$ of formal balls such that, for every $x \in X$, the set
  of formal balls $(y, s) \in \mathcal B$ such that
  $(y, s) \ll (x, 0)$ is directed, and has $(x, 0)$ as supremum.
%
\end{defi}

In the sequel, we fix a Yoneda-complete quasi-metric space $X, d$ with
a continuous poset of formal balls, and a local basis $\mathcal B$ of
$\mathbf B (X, d)$.
\begin{defi}
  \label{defn:BD}
  The poset $\mathbf B' (X, d; \mathcal B)$ is defined as follows.
  Its elements are all the formal balls of the form $(x, 0)$,
  $x \in X$, and those in the local basis $\mathcal B$.  Its ordering
  is defined by $(x, r) \sqsubseteq (y, s)$ if and only if
  $(x, r) \sqsubset (y,s)$ or $(x,r)=(y,s)$, where
  $(x, r) \sqsubset (y, s)$ if and only if:
  \begin{enumerate}
  \item either $(x, r) \ll (y, s)$ and $r \geq 2s$, where $\ll$ is the
    way-below relation in $\mathbf B (X, d)$,
  \item or $r=s=0$ and $x \leq^d y$.
  \end{enumerate}
\end{defi}
The second clause ensures that, equating $x \in X$ with
$(x, 0) \in \mathbf B' (X, d; \mathcal B)$, the ordering on $X$ is the
restriction of $\sqsubseteq$.  The first clause can be interpreted as
saying that to move up (strictly) among the elements of $\mathcal B$,
we must not only jump high---take an $(y, s)$ that is way-above
$(x, r)$---but also reduce radii by a constant factor.  We take $2$
for this factor, but this is arbitrary: any constant strictly larger
than $1$ would work equally well.

The first clause also allows one to compare an element $(x, r) \in
\mathcal B$ with an element of the form $(y, 0)$, not just to compare
two elements of $\mathcal B$.  However, one must note that $(x, r) \ll
(y, s)$ forces $r \neq 0$, hence $(x, r) \in \mathcal B$.  This is a
direct consequence of Lemma~\ref{lemma:+eps}:
\begin{rem}
  \label{rem:r!=0}
  If $(x, r) \ll (y, s)$ then $r > s$; in particular, $r \neq 0$.
\end{rem}

We shall say ``$\sqsubseteq$-directed'' or ``$\leq^{d^+}$-directed''
to make clear with respect to which ordering directedness is assumed,
and similarly for other epithets.
\begin{fact}
  \label{fact:dir}
  Plainly, $(x, r) \sqsubseteq (y, s)$ implies
  $(x, r) \leq^{d^+} (y, s)$, and that implies that every
  $\sqsubseteq$-directed family is also $\leq^{d^+}$-directed.
\end{fact}
The following technical lemma will be useful.
\begin{lem}
  \label{lemma:BD:trick}
  Let ${(x_i, r_i)}_{i \in I}$ be a $\sqsubseteq$-directed family in
  $\mathbf B' (X, d; \mathcal B)$, with $r_i \neq 0$ for every
  $i \in I$, and assume that it has no $\sqsubseteq$-largest element.
  For every $i \in I$, there is a $j \in I$ such that
  $(x_i, r_i) \ll (x_j, r_j)$ and $r_i \geq 2r_j$.
\end{lem}
\proof Since $(x_i, r_i)$ is not $\sqsubseteq$-largest, there is a
$(x_k, r_k)$ such that $(x_k, r_k) \not\sqsubseteq (x_i, r_i)$.  By
directedness, find $(x_j, r_j)$ such that
$(x_i, r_i), (x_k, r_k) \sqsubseteq (x_j, r_j)$.  It cannot be that
$(x_i, r_i) = (x_j, r_j)$, it cannot be either that $r_i=r_j=0$ and
$x_i \leq^d x_j$ since all radii are assumed non-zero, so
$(x_i, r_i) \ll (x_j, r_j)$ and $r_i \geq 2r_j$.  \qed

Equally useful is the following consequence, which should be
interpreted in the light of Example~\ref{exa:qideal:no}.  There is a
poset $B$ consisting of those elements of
$\mathbf B' (X, d; \mathcal B)$ whose radius is non-zero, ordered by
$\sqsubseteq$.  The lemma below states that the directed families $D$
in $B$ that have no largest element have no upper bound in $B$.  Hence
$\Idl (B)$ will be a quasi-ideal domain.  There is some remaining work
to do to show that $\Idl (B)$ is in fact isomorphic to
$\mathbf B' (X, d; \mathcal B)$, but this is a good start.
\begin{lem}
  \label{lemma:BD:trick:r=0}
  Let ${(x_i, r_i)}_{i \in I}$ be a $\sqsubseteq$-directed family in
  $\mathbf B' (X, d; \mathcal B)$, with $r_i \neq 0$ for every
  $i \in I$, and assume that it has no $\sqsubseteq$-largest element.
  Then $\inf_{i \in I} r_i = 0$.
\end{lem}
\proof Iterating Lemma~\ref{lemma:BD:trick} from some arbitrary index
$i_0 \in I$, we obtain a sequence
$(x_{i_0}, r_{i_0}) \ll (x_{i_1}, r_{i_1}) \ll \cdots \ll (x_{i_k},
r_{i_k}) \ll \cdots$
with $r_{i_0} \geq 2r_{i_1} \geq \cdots \geq 2^k r_{i_k} \geq \cdots$.
The infimum of those values is $0$.  \qed

\begin{lem}
  \label{lemma:BD:dcpo}
  $\mathbf B' (X, d; \mathcal B)$ is a dcpo, and directed suprema are
  computed exactly as in the larger dcpo $\mathbf B (X, d)$.
\end{lem}
\proof Let ${(x_i, r_i)}_{i \in I}$ be a $\sqsubseteq$-directed family
in $\mathbf B' (X, d; \mathcal B)$, and let $(x, r)$ be its
$\leq^{d^+}$-supremum, in $\mathbf B (X, d)$.  If that supremum is
reached, namely if $(x, r) = (x_i, r_i)$ for some $i \in I$, then we
claim that $(x, r)$ is also the $\sqsubseteq$-supremum.  To show that
it is an $\sqsubseteq$-upper bound, we must show that
$(x_j, r_j) \sqsubseteq (x, r)$, for any $j \in I$.  By directedness,
find $k \in I$ so that
$(x_i, r_i), (x_j, r_j) \sqsubseteq (x_k, r_k)$.  Then
$(x_i, r_i) \leq^{d^+} (x_k, r_k) \leq^{d^+} (x, r) = (x_i, r_i)$, so
$(x_i, r_i) = (x_k, r_k)$.  That implies
$(x_j, r_j) \sqsubseteq (x_i, r_i) = (x, r)$.  It is the
$\sqsubseteq$-supremum because it is attained.


If $r_j=0$ for some $j \in I$, then the subfamily of those formal
balls $(x_i, r_i)$ such that $(x_j, r_j) \sqsubseteq (x_i, r_i)$ is
again $\sqsubseteq$-directed, and has the same $\sqsubseteq$-upper
bounds and the same $\leq^{d^+}$-upper bounds as the original family.
If $r_j=0$, $(x_j, r_j) \sqsubseteq (x_i, r_i)$ implies $r_i=0$, so
that subfamily consists of formal balls of radius $0$.  On formal
balls of radius $0$, $\sqsubseteq$ and $\leq^{d^+}$ coincide, because
of Remark~\ref{rem:r!=0}.  Therefore $(x, r)$ is the
$\sqsubseteq$-supremum of our original family again.  Note that, since
$r=0$, it is, in particular, in $\mathbf B' (X, d; \mathcal B)$.

Henceforth assume that the family has no $\sqsubseteq$-largest
element, and that $r_j \neq 0$ for every $j \in I$.  In particular,
the whole family is inside $\mathcal B$.  For every $i \in I$, use
Lemma~\ref{lemma:BD:trick} and find $j \in I$ such that
$(x_i, r_i) \ll (x_j, r_j)$ and $r_i \geq 2r_j$.  Since
$(x_j, r_j) \leq^{d^+} (x, r)$, it follows that
$(x_i, r_i) \ll (x, r)$ and $r_i \geq 2r$, whence
$(x_i, r_i) \sqsubset (x, r)$, at least if we can show that $(x, r)$
is in $\mathcal B' (X, d; \mathcal B)$, which we shall do below.

Let $(y, s)$ be an arbitrary $\leq^{d^+}$-upper bound of the family.
Note that $s \leq \inf_{i \in I} r_i$, and that
$\inf_{i \in I} r_i = 0$ by Lemma~\ref{lemma:BD:trick:r=0}.  This
applies to $(y, s)=(x,r)$, so $r=0$ as well.  That shows that $(x, r)$
is in $\mathcal B' (X, d; \mathcal B)$, as promised, so $(x, r)$ is a
$\sqsubseteq$-upper bound of the family.

Now consider any $\sqsubseteq$-upper bound $(y, s)$ of the family.
This is, in particular, a $\leq^{d^+}$-upper bound, so $s=0$, as we
have just seen.  Since $(x, r)$ is the least $\leq^{d^+}$-upper bound,
$(x, r) \leq^{d^+} (y, s)$.  We now use $r=s=0$ to conclude that
$(x, r) \sqsubseteq (y, s)$.  \qed

\begin{lem}
  \label{lemma:approx}
  Every element $(x, 0)$, $x \in X$, is the supremum in $\mathbf B'
  (X, d; \mathcal B)$ of a $\sqsubseteq$-directed family of elements
  $(y, s) \in \mathcal B$ with $(y, s) \ll (x, 0)$.
\end{lem}
\proof Since $\mathcal B$ is a local basis of $\mathbf B (X, d)$,
$(x, 0)$ is the $\leq^{d^+}$-supremum of a $\leq^{d^+}$-directed
family $D = {(x_i, r_i)}_{i \in I}$ of formal balls in $\mathcal B$,
with $(x_i, r_i) \ll (x, 0)$ for every $i \in I$.  In particular,
$(x_i, r_i) \sqsubseteq (x, 0)$.  There is no reason why $D$ should be
$\sqsubseteq$-directed, but we repair this by a variant of a standard
trick, allowing us to obtain a subfamily of $D$ that is
$\sqsubseteq$-directed and will have the same supremum.

Note that $r_i \neq 0$ for each $i \in I$, by Remark~\ref{rem:r!=0}.
Also, since $(x, 0)$ is the $\leq^{d^+}$-supremum of $D$,
$0 = \inf_{i \in I} r_i$.  (Recall that $X, d$ is Yoneda-complete,
hence standard.)

Our new family consists of formal balls which we write $(x_J, r_J)$,
indexed by finite subsets $J$ of $I$, and is constructed by induction
on the cardinality of $J$.  We will ensure that each $(x_J, r_J)$ is
in $D$, and that $(x_J, r_J) \ll (x, 0)$, for every $J$.  We shall
also make sure that $(x_K, r_K) \sqsubseteq (x_J, r_J)$ for every
$K \subseteq J$, and that will ensure that our new family is
$\sqsubseteq$-directed.

When $J=\emptyset$, we pick $(x_J, r_J)$ arbitrarily from $D$.
Otherwise, consider $U$, the intersection of the subsets
$\uuarrow (y, s)$, $(y, s) \in J$, and of the subsets of the form
$\uuarrow (x_K, r_K)$ for $K \subsetneq J$.  All those subsets are
Scott-open, because $\mathbf B (X, d)$ is a continuous poset.  As a
finite intersection of open subsets, $U$ is open.  Since
$(x_i, r_i) \ll (x, 0)$ for every $i \in I$, and by induction
hypothesis, $U$ contains $(x, 0)$.  Hence some $(x_i, r_i)$ is in $U$.
Since $\inf_{j \in I} r_j = 0$ and $r_i \neq 0$, we can also find
$(x_j, r_j)$ in $D$ so that $r_i \geq 2 r_j$, and by
$\leq^{d^+}$-directedness, we can even require that
$(x_i, r_i) \leq^{d^+} (x_j, r_j)$.  Now decide to set
$(x_J, r_J) = (x_j, r_j)$.  (All that, naturally, uses the Axiom of
Choice rather heavily.)  That is enough to ensure that
$(x_K, r_K) \sqsubset (x_J, r_J)$ for every $K \subsetneq J$, since
$(x_K, r_K) \ll (x_i, r_i) \leq^{d^+} (x_j, r_j)$ and
$r_K \geq r_i \geq 2r_j$.

By construction, $(x_i, r_i) \sqsubseteq (x_{\{i\}}, r_{\{i\}})$, too,
so our new family of formal balls $(x_J, r_J)$ is not only
$\sqsubseteq$-directed, but also $\sqsubseteq$-cofinal in $D$, hence
it has $(x, 0)$ as $\sqsubseteq$-supremum.  \qed

\begin{lem}
  \label{lemma:finite}
  The finite elements of $\mathbf B' (X, d; \mathcal B)$ are its
  elements with non-zero radius.  (In particular, they are all in
  $\mathcal B$.)
\end{lem}
\proof Let $(x, r) \in \mathcal B$, with $r \neq 0$.  We show that
$(x, r)$ is $\sqsubseteq$-finite.  Consider a $\sqsubseteq$-directed
family ${(x_i, r_i)}_{i \in I}$ with supremum $(y, s)$ such that
$(x, r) \sqsubseteq (y, s)$.  If the supremum is attained, say at
$i \in I$, then $(x, r) \sqsubseteq (x_i, r_i)$, and we are done.

Otherwise, we note that since $r \neq 0$ and
$(x, r) \sqsubseteq (y, s)$, either $(x, r) = (y, s)$ or
$(x, r) \ll (y, s)$ and $r \geq 2s$.  In the latter case, for some
$i \in I$, and using interpolation, $(x, r) \ll (x_i, r_i)$.

We now, again, enquire whether $r_i=0$ for some $i$.
If that is the case, then we can without loss of generality assume
that $r_j=0$ for every $j \in I$, taking a cofinal subfamily if
necessary.  Note that $s = \inf_j r_j = 0$.
Since $(x, r) \sqsubseteq (y, s)$ and $r \neq 0$, a case analysis on
the definition of $\sqsubseteq$ shows that $(x, r) \ll (y, s)$.
Therefore, and using interpolation, $(x, r) \ll (x_i, r_i)$ for some
$i \in I$.  Since $r_i=0$, $(x, r) \sqsubseteq (x_i, r_i)$ trivially.

If the supremum is not attained and $r_i \neq 0$ for every $i \in I$,
then $s \leq \inf_{i \in I} r_i = 0$, using
Lemma~\ref{lemma:BD:trick:r=0}.  Since $r \neq 0$ and $s=0$, as above
we find an index $i \in I$ such that $(x, r) \ll (x_i, r_i)$.  Now use
Lemma~\ref{lemma:BD:trick}: there is an index $j \in I$ such that
$(x_i, r_i) \ll (x_j, r_j)$ and $r_i \geq 2r_j$.  It follows that
$(x, r) \ll (x_j, r_j)$ and $r \geq 2r_j$, namely,
$(x, r) \sqsubset (x_j, r_j)$.

Next, we show that no element of the form $(x, 0)$ is
$\sqsubseteq$-finite.  Use Lemma~\ref{lemma:approx} to find a
$\sqsubseteq$-directed family ${(x_i, r_i)}_{i \in I}$ whose supremum
is $(x, 0)$, and with $(x_i, r_i) \ll (x, 0)$.  By
Remark~\ref{rem:r!=0}, $r_i \neq 0$.  If $(x, 0)$ were
$\sqsubseteq$-finite, there would be an $i \in I$ such that
$(x, 0) \sqsubseteq (x_i, r_i)$.  In particular,
$(x, 0) \leq^{d^+} (x_i, r_i)$, which implies $r_i=0$: that is
impossible.  \qed

\begin{prop}
  \label{prop:qideal}
  $\mathbf B' (X, d; \mathcal B)$ is a quasi-ideal domain.
\end{prop}
\proof By Lemma~\ref{lemma:approx} and Lemma~\ref{lemma:finite}, it is
algebraic.  Also, if $(x, r)$ is below a finite element $(y, s)$,
namely one with $s \neq 0$, then $r \geq s$, which implies that
$r\neq 0$, hence that $(x, r)$ is finite, too.  \qed

\begin{prop}
  \label{prop:qideal:model}
  Let $X, d$ be a continuous Yoneda-complete quasi-metric space, and
  let $\mathcal B$ be a local basis of $\mathbf B (X, d)$.  Then $X$
  with the $d$-Scott topology is homeomorphic to the subspace of limit
  elements of the quasi-ideal domain $\mathbf B' (X, d; \mathcal B)$.
\end{prop}
\proof Temporarily call \emph{$d$-Scott'} topology the topology
induced by the Scott topology on $\mathbf B' (X, d; \mathcal B)$ on
its subspace $X$.  In the light of Proposition~\ref{prop:qideal}, it
only remains to show that the $d$-Scott and the $d$-Scott' topologies
are the same.

Let $U$ be a $d$-Scott open in $X$, and write it as $V \cap X$, where
$V$ is Scott-open in $\mathbf B (X, d)$.  Let
$V' = V \cap \mathbf B' (X, d; \mathcal B)$.  Clearly,
$V' \cap X = U$.  We claim that $V'$ is Scott-open in
$\mathbf B' (X, d; \mathcal B)$.  It is $\sqsubseteq$-upwards-closed
because $V$ is $\leq^{d^+}$-upwards-closed, and since $\sqsubseteq$
implies $\leq^{d^+}$ (Fact~\ref{fact:dir}).  For a
$\sqsubseteq$-directed family $D$ in $\mathbf B' (X, d; \mathcal B)$
with $\sqsubseteq$-supremum $(x, r)$ in $V'$, $D$ is
$\leq^{d^+}$-directed and has $(x, r)$ as $\leq^{d^+}$-supremum, by
Lemma~\ref{lemma:BD:dcpo}.  Since $(x, r)$ is in $V$, some element of
$D$ is in $V$, and therefore necessarily in $V'$ as well.  Since $V'$
is Scott-open in $\mathbf B' (X, d; \mathcal B)$, $U = V' \cap X$ is
then $d$-Scott' open.

Conversely, let $U$ be $d$-Scott' open, and let $V'$ be a Scott-open
subset of $\mathbf B' (X, d; \mathcal B)$ such that $U = V' \cap X$.
For every $x \in U$, $(x, 0)$ is in $V'$, therefore $(x, 0)$ has an
open neighborhood $\upc (y, r)$ included in $V'$, where $(y, r)$ is a
finite element of $\mathbf B' (X, d; \mathcal B)$, and $\upc$ denotes
upwards-closure in $\mathbf B' (X, d; \mathcal B)$.  Since $(y, r)$ is
finite, $r$ is non-zero (Lemma~\ref{lemma:finite}).  In particular,
since $(y, r) \sqsubseteq (x, 0)$, we must have $(y, r) \ll (x, 0)$ in
$\mathbf B (X, d)$.  Therefore $\uuarrow (y, r) \cap X$ is a $d$-Scott
open neighborhood of $x$.

For every $z \in \uuarrow (y, r) \cap X$, we have $(y, r) \ll (z, 0)$,
so $(y, r) \sqsubseteq (z, 0)$.  Therefore,
$(z, 0) \subseteq \upc (y, r) \subseteq V'$, and that entails
$z \in U$.  Hence $\uuarrow (y, r) \cap X$ is included in $U$.
Since $U$ is a $d$-Scott open neighborhood of each of its points, $U$
is $d$-Scott open.  \qed

Taking $\mathcal B = \mathbf B (X, d)$ itself, for example, we obtain:
\begin{thm}
  \label{thm:qideal:model}
  Every continuous Yoneda-complete quasi-metric space, in its
  $d$-Scott topology, has a quasi-ideal model.
\end{thm}

Looking back at Example~\ref{ex:Scott=dScott}, we obtain the following.
\begin{cor}
  \label{corl:qideal:model:cpo}
  Every continuous dcpo has a quasi-ideal model.
\end{cor}

Making some countability assumptions allows us to refine those
results.  Recall that a dcpo is $\omega$-continuous if and only if it
is continuous and has a countable basis.  Equivalently, a continuous
dcpo is $\omega$-continuous if and only if its Scott topology is
countably based \cite[Proposition~3.1]{Norberg:randomsets}.

Applying Proposition~\ref{prop:qideal:model} to a countable basis
$\mathcal B$ of $\mathbf B (X, d)$, $\mathbf B' (X, d; \mathcal B)$
will be a quasi-ideal domain with countably many finite elements, that
is, an $\omega$-quasi-ideal domain.  Hence:
\begin{prop}
  \label{prop:qideal:omegamodel}
  Every Yoneda-complete quasi-metric space with an $\omega$-continuous
  poset of formal balls, equipped with its $d$-Scott topology, has an
  $\omega$-quasi-ideal model.

  In particular, every $\omega$-continuous dcpo has an
  $\omega$-quasi-ideal model.
\end{prop}
Matthew de Brecht introduced the quasi-Polish spaces as the
topological spaces underlying the countably based Smyth-complete
spaces in their open ball topology \cite{deBrecht:qPolish}
(equivalently, in their $d$-Scott topology, since the two topologies
coincide on Smyth-complete spaces).  Among the many equivalent
characterizations of quasi-Polish spaces, Theorem~53 of op.\ cit.\
states that the quasi-Polish spaces are also the topological spaces
that embed into an $\omega$-algebraic dcpo as its subset of non-finite
elements.  The equivalence between (1) and (4) below strengthens this
result.

\begin{thm}
  \label{thm:qPolish}
  The following are equivalent for a topological space $X$:
  \begin{enumerate}
  \item $X$ is quasi-Polish;
  \item $X$ is countably based and there is a quasi-metric $d$ that
    makes $X, d$ a Smyth-complete quasi-metric space;
  \item there is a quasi-metric $d$ that makes $X, d$ a
    Yoneda-complete quasi-metric space with an $\omega$-continuous
    poset of formal balls;
  \item $X$ has an $\omega$-quasi-ideal model.
  \end{enumerate}
\end{thm}
\proof The previous discussion shows the implication $(4) \limp (1)$,
in particular.  $(1) \limp (2)$ is by definition.  $(3) \limp (4)$ is
by Proposition~\ref{prop:qideal:omegamodel}.  It remains to show
$(2) \limp (3)$.  Assume $X, d$ Smyth-complete.  In particular, it is
Yoneda-complete and continuous.  We need to show that
$\mathbf B (X, d)$ is $\omega$-continuous, and we exhibit a countable
basis to that end.

Let $d^{sym}$ be the symmetrization of $d$:
$d^{sym} (x, y) = \max (d (x, y), d (y, x))$.  Since $X$ is countably
based in the open ball topology of $d$ (which coincides with the
$d$-Scott topology), it has a countable subset $D$ that is dense with
respect to the open ball topology of $d^{sym}$.  In fact, this is an
equivalence, due to H.-P. A. K\"unzi \cite{Kunzi:quasiuniform}, and
mentioned as Proposition~13 in \cite{deBrecht:qPolish} (see also
\cite[Lemma~6.3.48]{JGL-topology}): a quasi-metric space has a
countably-based open ball topology if and only if its symmetrization
has a countable dense subset.  Let $\mathcal B$ be the countable set
of formal balls $(x, r)$ with $x \in D$ and $r$ rational.

Fix a formal ball $(y, s)$, and consider
$A = \{(x, r) \in \mathcal B \mid d^{sym} (x, y) < r-s\}$.  Given
finitely many elements $(x_1, r_1)$, \ldots, $(x_n, r_n)$ of $A$, let
$\epsilon$ be any positive number less than or equal to
$\min_{i=1}^n (r_i - s - d^{sym} (x_i, y))$ such that $r = s+\epsilon/2$
is rational.  By density, $B^{d^{sym}}_{y, <\epsilon/2}$ contains some
$x \in D$.  By symmetry, $d^{sym} (x, y) < \epsilon/2 = r-s$, so
$(x, r)$ is in $A$.  Moreover, for each $i$,
$d^{sym} (x_i, x) \leq d^{sym} (x_i, y) + d^{sym} (y, x) < d^{sym}
(x_i, y) + \epsilon/2 \leq (r_i - s - \epsilon) + \epsilon/2 = r_i - s
- \epsilon/2 = r_i - s - (r-s) = r_i-r$.
In particular, $d (x_i, x) \leq r_i - r$, which entails
$(x_i, x) \leq^{d^+} (x, r)$.  It follows that $A$ is a directed
family of elements.

Every element $(x, r)$ of $A$ is such that $d (x, y) < r-s$, namely
$(x, r) \prec (y, s)$.  By the Romaguera-Valero Theorem,
$\ll = \prec$, so every element of $A$ is way-below $(y, s)$.  The net
${(x)}_{(x, r) \in A, \leq^{d^+}}$ is Cauchy-weightable by definition,
and $y$ is its limit in the open ball topology of $X, d^{sym}$: indeed
every open neighborhood of $y$ in that topology contains an open ball
$B^{d^{sym}}_{y, <\eta}$ for some $\eta > 0$, which itself contains
some $x \in D$ by density; $A$ must contain $(x, r)$ for some rational
positive $r$, namely with $r$ such that $0 < r-s \leq \eta$.  Every
limit in the open ball topology of $d^{sym}$ is a $d^{sym}$-limit
\cite[Proposition~7.1.19]{JGL-topology}.  It follows that the supremum
of $A$ is the formal ball whose center is the $d^{sym}$-limit of
${(x)}_{(x, r) \in A, \leq^{d^+}}$, which is $y$ as we have just seen,
and whose radius is $\inf \{r \mid (x, r) \in A\} = s$.  Hence
$\sup A = (y, s)$.

Since every formal ball $(y, s)$ is the supremum of a directed family
$A$ of elements from the countable set $\mathcal B$ and way-below
$(y, s)$, $\mathcal B$ is a countable basis of $\mathbf B (X, d)$.
\qed


\section{Conclusion and Open Problems}
\label{sec:concl-open-probl}

We have shown a variety of results on quasi-metric spaces, and all
share one feature: they are all proved domain-theoretically, by
reasoning on the poset of formal balls.  This proves to be a useful
complement to the view of quasi-metric spaces as enriched categories,
and works by relatively simple reductions to notions and techniques
from ordinary domain theory.

Some questions remain open, as usual.  Is there any form of converse
to Theorem~\ref{thm:qideal:model}?  In general, what are the spaces
that have a quasi-ideal model?  Theorem~\ref{thm:qPolish} answers the
question completely for countably based spaces, but what about
non-countably based spaces?  Keye Martin showed that the metric spaces
that have an ideal model are exactly the complete metric spaces
\cite{Martin:idealmodels}.  However, ideal models are not only
algebraic, but also first-countable, and that is crucial.  There is no
reason to believe that quasi-ideal models are first-countable, and
continuous Yoneda-complete quasi-metric spaces are not in general
first-countable in their $d$-Scott topology either (as they contain
all continuous dcpos already, see Example~\ref{ex:<=} and
Example~\ref{ex:Scott=dScott}).



\bibliographystyle{alpha}

\newcommand{\etalchar}[1]{$^{#1}$}



\end{document}